\documentclass[10pt,a4paper]{article}

\usepackage{amssymb,amsmath,amsthm,amscd,enumerate}
\usepackage[all]{xy}
\usepackage[french]{babel}
\newcommand{\C}{{\bf C}}
\newcommand{\R}{{\bf R}}
\newcommand{\U}{{\bf U}}
\newcommand{\p}{{\bf M}}
\newcommand{\Z}{{\bf Z}}
\newcommand{\pp}{\mathfrak{p}}
\newcommand{\Sp}{{\rm Spec}}
\newcommand{\Q}{{\bf Q}}

\advance\oddsidemargin -1,5cm
\advance\evensidemargin -2cm

\begin{document}
\begin{center}
{\Large{\bf Un module inversible associ\'e au ruban de M\"obius, 

et quelques autres}}
\vspace{0,4cm}

Daniel Ferrand

{\it mars 2007}

\end{center}

\vspace{3cm}

{\bf Abstract}\; Invertible modules are, in a sense, the simplest interesting modules one encounters in commutative algebra, and one cannot avoid them when dealing with algebraic problems coming from number theory or geometry; they may be seen as the algebraic counterpart of the notion of twisting, or  glueing in geometry.

Unfortunately, in most of the textbooks, invertible modules are introduced after quite an amount of prerequisites which dissuade students from discovering  them, if they don't intend  to become specialists in algebraic geometry; this notion deserves,  however, to be part of the basic education  of any mathematician.

To remedy this fact a little, this paper, and the lectures it is based on, aims to introduce the notion of invertible module as directly as possible. The text begins, in a very elementary way, by attaching such a module to the famous M\"obius strip --- mere translation exercises between a geometric language dealing with a "concrete" object, and perhaps more abstract algebraic notions.
In the paragraph 4, invertible modules are defined as projective modules of rank one; those generated by two elements are constructed in an explicit (versal) way. In the next paragraph some  rudiments of Zariski topology are required  in order to identify an invertible A-module as a module  which is locally isomorphic to A. 

In the paragraph 6, I give a new way of attaching to a binary quadratic form an invertible module over a suitable ring of integers; such a correspondence between forms and \emph{ideals} was already known by Gauss and Dedekind, but it is drastically simplified by the use of the notion of invertible module which is more flexible than that of ideal.

Then follows a very short proof of the triviality of invertible modules over factorial rings.

After that, in the paragraphs 8 and 9,  the level of the text is raising a little : tensor products are taken for granted and I review, and use, the first steps of descent theory.

The paper ends with specialisation to Galois coverings, and invertible modules constructed from cocycles  of the Galois group ; the particular case where the group is abelian receives some attention in the paragraph 11.

None of the results is really new, but the short ways I manage to get them deserve, I think,  to be known. The first half of the text may be understood by students at the advanced undergraduate level.
\newpage

\noindent Introduction
\medskip

\noindent 1. \,G\' eom\' etrie : le ruban de M\"obius
\medskip

\noindent 2. \,Traduction en termes alg\' ebriques
\medskip

\noindent 3. \,Alg\`ebre : l'anneau $A = \R[X, Y]/(X^2+Y^2-1)$
\medskip

\noindent 4. \,Modules inversibles

\medskip

\noindent 5. \,Modules localement isomorphes \`a l'anneau 
\medskip

\noindent 6. \,Application : le module inversible associ\'e \`a une forme quadratique binaire
\medskip

\noindent 7. \,Trivialit\'e des modules inversibles sur un anneau factoriel
\medskip

\noindent 8. \,Modules inversibles et produit tensoriel
\medskip

\noindent 9. \,Descente
\medskip

\noindent 10. Constructions galoisiennes
\medskip

\noindent 11. Application : la th\'eorie de Kummer

\vspace{0,5cm}

\begin{center}
{\bf Introduction}
\end{center}
\vspace{0,5cm}

Ce texte est centr\'e sur la notion de \emph{module inversible} ; il en explore la plupart des aspects alg\'ebriques, mais de fa\c{c}on \'el\'ementaire. Si aucun r\'esultat n'est vraiment nouveau, les parcours - directs - pour y parvenir le sont peut-\^etre.

J'ai r\'edig\'e ce texte, qui prolonge des cours de Master anciens, en pensant \`a des lecteurs qui n'approfondiront peut-\^etre pas la g\'eom\'etrie alg\'ebrique, et pour lesquels on peut (on doit) donc \'eviter les longs pr\'eliminaires dissuasifs des pr\'esentations usuelles.
\medskip 

Les trois premiers paragraphes sont tr\`es \'el\'ementaires : ils n'utilisent m\^eme pas  le produit tensoriel, qui ne sera employ\'e qu'\`a partir du paragraphe 8. Ce d\'ebut  est une sorte de variation g\'eom\'etrique et alg\'ebrique sur le th\`eme du ruban de M\"obius ; j'en donne un plongement dans $\C^2$ qui \'evite l'usage habituel des fonctions trigonom\'etriques, et qui est plus commode que la d\'efinition comme quotient de $\R^2$ par l'action d'un groupe.

Le paragraphe 4 d\'efinit un module inversible comme un module projectif de rang 1 ; d\'efinition et r\'esultats utilisent uniquement  l'alg\`ebre lin\'eaire sur un anneau. Les modules inversibles engendr\'es par deux \'el\'ements sont construits tr\`es explicitement.

Le paragraphe 5 emploie quelques mots provenant de la topologie (de Zariski) pour d\'esigner des constructions purement alg\'ebriques qui auraient eu plus de sens si le dictionnaire entre elles et cette topologie avait \'et\'e developp\'e comme il le m\'erite. Mais cela aurait contrecar\'e  le projet de ce texte : faire court.

Le paragraphe 6 sugg\`ere que le concept de module inversible peut simplifier drastiquement la pr\'esentation v\'en\'erable des liens entre formes quadratiques binaires et classes d'id\'eaux d'entiers alg\'ebriques. L'id\'ee est probablement nouvelle.

Au paragraphe 7, je d\'emontre en quelques lignes le r\'esultat, frappant, de son titre.

Les quatre derniers paragraphes, peut-\^etre un peu plus \'elabor\'es, utilisent le produit tensoriel et les diagrammes commutatifs. J'ai tenu \`a mettre en \oe{}uvre (une petite partie de)  la th\'eorie de la descente pour en montrer l'efficacit\'e, et parce que je regrette qu'elle n'ait pas encore acquis sa place dans l'enseignement, parmi les outils de base. Le paragraphe sur les constructions galoisiennes s'ach\`eve sur le constat que la construction du d\'ebut du texte en est une ! 

Le cas o\`u le groupe de Galois est ab\'elien est pr\'ecis\'e dans le dernier paragraphe.
\vspace{2cm}

{\Large {\bf 1. G\' eom\' etrie : le ruban de M\"obius.}}\\

Sur votre bande de papier calque, de forme rectangulaire assez along\' ee,
tracez une
droite parall\`ele aux deux grands c\^ot\' es, et \' equidistante d'iceux ; en chaque
point de
cette droite pensez \`a, ou tracez, un segment perpendiculaire, que vous
imaginerez de
longueur infinie. Collez les petits c\^ot\' es du rectangle, apr\`es avoir effectu\' e \emph{le
demi-tour fatidique}. La droite parall\`ele aux grands c\^ot\'es devient un cercle, et les deux grands c\^ot\'es sont rabout\'es et forment l'unique bord de cette surface. Cette chose, la
bande de
M\"obius, devrait vous faire penser \`a une famille de droites param\'
etr\' ee par le cercle, famille radicalement diff\' erente de celle repr\'esent\'ee par un cylindre... et
pourtant il
suffit de les couper le long d'une droite
verticale
pour que ces deux surfaces deviennent \' evidemment isomorphes.\\

Cet objet g\' eom\' etrique a un analogue alg\' ebrique. Consid\' erez l'anneau $A = \R[X,
Y]/(X^2+Y^2-1)$  des applications polynomiales d\' efinies
sur le cercle, \`a valeurs r\' eelles. C'est un anneau int\' egralement clos mais l'id\'eal engendr\'e par $1+x$ et $y$ \emph{n'est pas principal}, et ceci bien que le complexifi\'e de $A$, $\C[X,
Y]/(X^2+Y^2-1)$, soit un anneau principal.\\

Dans le paragraphe 1, je donne un mod\`ele math\' ematique tr\`es simple du
ruban de
M\"obius. Certaines propri\'et\'es alg\'ebriques de ce mod\`ele conduisent, au \S 2,  \`a lui associer un module sur l'anneau $A$, module dont le ruban peut \^etre vu comme la figure g\'eom\'etrique ; le fait que le ruban ne soit pas un cylindre est traduit par le fait que ce module n'est pas libre. Au \S 3 on \'etablit un isomorphisme entre cet  $A$-module et l'id\'eal de $A$ engendre par $1+x$  et $y$.\bigskip

L'article de {\sc Swan} $[8]$ contient, p. 273, quelques lignes qui m'ont \'eclair\'e.

\vspace{0,5cm}

Dans tout ce texte, le symbole
{\bf U}  d\' esigne le groupe des nombres complexes de module $1$, vu,
g\' eom\' etriquement, comme le cercle unit\' e.\\

{\bf 1.1. D\' efinition}\\

Le ruban de M\"obius est le sous-ensemble de $
\U
\times \C$ d\'efini par :
$$
\p \; = \; \{ (u, z) \in \U \times \C , u\bar{z} = z\} .
$$
On le munit de la topologie induite par celle de  $\U\times \C$.
On utilisera constamment la\\

{\bf Remarque-cl\' e 1.1.1.}\; {\it Si $v \in \U$, alors la relation
$v^2\bar{z} =
z$  \' equivaut \`a
$ z = \mu v$, avec $\mu \in \R$}.
\bigskip

\noindent En effet, la relation en question s'\' ecrit aussi $ \overline{(z/v)} =
z/v$.$\Box$
\vspace{0,5cm}

\noindent La premi\`ere projection $ \U \times \C  \longrightarrow \U$ induit une
application continue
$$
p : \p \quad \longrightarrow \quad \U ,\qquad (u, z) \mapsto u .
$$
Pour tout $u \in \U$, la fibre $\p_u = p^{-1}(u)$ est un $\R$-espace
vectoriel de dimension $1$. En effet, \C \, \' etant alg\' ebriquement clos, tout nombre 
complexe admet une racine carr\' ee, et si on choisit une racine carr\' ee de
$u$, soit un \' el\' ement $ v\in \U$ tel que $v^2 = u$, alors, d'apr\`es 1.1.1,
on a
$$
\p_u \; = \; \R v  \subset \C.
$$
La droite $\p_u$ est donc la bissectrice de l'angle $\widehat{(1, u)}$ .
Ainsi, l'application $p : \p \longrightarrow \U$ conduit \`a voir $\p$
comme une famille de droites vectorielles dans le plan, param\' etr\' ees par le
cercle unit\' e, et qui \, \og tournent deux fois moins vite\fg \, que le
param\`etre.\\

Avant de faire le lien entre cette d\' efinition et la bande de M\"obius
\og concr\`ete\fg, il faut indiquer comment on obtient le ruban de
M\"obius comme quotient du cylindre $\U \times \R$ par l'action d'un
groupe d'ordre $2$.\\

{\bf Lemme 1.1.2.}\;{\it Le carr\' e suivant est cart\' esien, i.e il permet
d'identifier} $\U\times \R$ {\it au sous-ensemble de $\U \times \p$ form\' e
des couples} $(v, m)$ {\it tels que} $v^2 = p(m)$.
$$\begin{CD}
\U\times \R @>{(v, \mu)\mapsto (v^2, \mu v)}>> \p \\
@V{\text pr_1}VV  @VV p V\\
\U @>>{v\mapsto v^2}> \U
\end{CD}
$$
{\it De plus, l'application indiqu\' ee}\, $\U\times \R \rightarrow \p,\; (v,
\mu)\mapsto (v^2,
\mu v)$\, {\it permet d'identifier}\, $\p$\, {\it avec le quotient de}\,
$\U\times \R$ {\it sous l'involution} \, $(v, \mu) \mapsto (-v, -\mu)$.\\

(En termes savants  on
\' enoncerait que le ruban de M\"obius est un fibr\' e en droites sur
$\U$, trivialis\' e par le rev\^etement $v \mapsto v^2$ ; \og trivialis\'
e\fg\ veut dire : rendu isomorphe au cylindre $\U\times \R$.)

La premi\`ere assertion du lemme r\' ep\`ete la remarque-cl\' e. Pour v\' erifier la
seconde,
consid\' erons des \' el\' ements $(v, \mu)$  et $(v', \mu')$ de $\U\times \R$
ayant la m\^eme
image dans $\p$ ; on a donc $v^2 = {v'}^2$, d'o\`u $v' = \varepsilon v$, avec
$\varepsilon = \pm 1$, et
$\mu v = \mu' v' = \mu' \varepsilon v$ ; d'o\`u $(v', \mu') = \varepsilon
(v, \mu)$.
\vspace{0,5cm}

{\bf 1.2. Reconna\^itre la bande de M\"obius}\\

Rappelons la d\' efinition classique du ruban de M\"obius comme quotient
du plan $\R^2$ par l'action d'un groupe $\Gamma$ ({\sc Godbillon}, p.43 ;
{\sc Stillwell} p.25, o\`u la {\it bande} de M\"obius et nomm\' ee "M\"obius
strip", tandis que son {\it ruban} est nomm\' e "twisted cylinder", ce
qui est tr\`es \' evocateur).\\

 On consid\`ere la \, \og sym\' etrie gliss\' ee \fg \, ("glide
reflection")
\, $\gamma : \R^2 \rightarrow \R^2$\, d\' efinie par $\gamma(\lambda, \mu) =
(\lambda +1,
-\mu)$, et le groupe $\Gamma$ (isomorphe \`a {\bf Z}) engendr\' e par $\gamma$.

On remarque que $\gamma^2(\lambda, \mu) = (\lambda + 2, \mu)$ ; il est
donc naturel de passer d'abord au quotient par le sous-groupe engendr\' e
par
$\gamma^2$ puisqu'il n'agit que par translation sur le premier facteur
; l'application exponentielle donne alors un hom\' eomorphisme
$$
\R^2/<\gamma^2> \quad \widetilde{\longrightarrow}\quad \U\times \R, \qquad
(\lambda, \mu) \mapsto (e^{i\pi \lambda}, \mu).
$$
Comme $e^{i \pi} = -1$, la sym\' etrie gliss\' ee $\gamma$ induit sur le
quotient l'application $(u, \mu) \mapsto (-u, -\mu)$ ; c'est pr\' ecisement
l'involution consid\' er\' ee dans le lemme pr\' ec\' edent, et ce lemme implique donc
qu'on a un hom\' eomorphisme
$$
\R^2 / \Gamma \quad \widetilde{\longrightarrow}\quad \p.
$$

Pour remonter du \emph{ruban} (surface non compacte, sans bord, non
immergeable dans
l'espace) \`a la \emph{bande} (surface compacte, \`a bord, mais qu'on peut
fabriquer et
voir), on se restreint au rectangle
$D = [0, 1]\times [-1, 1]
\subset
\R^2$ ; les seuls points de $D$ que $\gamma$ envoie dans $D$ sont les
points du bord $(0, \mu)$, et $\gamma (0, \mu) = (1, -\mu)$ ; en
identifiant pour chaque $\mu \in [-1, 1]$, les points $(0, \mu)$  et  $(1,
-\mu)$ on
trouve la fameuse bande.

\newpage

{\bf 1.3. Les sections}\\

Une {\it section}  de $p$ est, par d\' efinition, une application continue
$s : \U \rightarrow \p$ telle que $p\circ s = {\rm Id}_{\U}$.
Il est habituel
d'\' ecrire ces applications verticalement :
$$
\xymatrix{
{\p}\ar@<1ex>[d]^p\\
{\U}\ar@<1ex>[u]^s
}
$$

\noindent La \emph{section nulle} est l'application $s_0 : \U \rightarrow
\p$ d\' efinie par
$s_0(u) = (u, 0)$.
\medskip

Une section peut \^etre vue comme la donn\' ee, pour chaque $u$,
d'un point
sur la droite $\p_u$, point qui varie contin\^ument avec $u$.\\

 Soit $s$ une  section ; notons $ f : \U \longrightarrow \C$ l'application compos\'ee
$$
\xymatrix{
{\p}\ar[r] &\U\times \C \ar[r]^{{\rm pr_2}} &\C \\
{\U}\ar[u]^s \ar[rru]_f
}
$$
L'application  $f$ est continue, tout comme $s$, et la condition $p\circ s = {\rm Id}_{\U}$ montre que pour tout $u \in \U$, on a 
$$s(u) = (u, f(u)).$$

Comme $(u, f(u))$ est dans $\p$, on a 
$$
u\overline{f(u)} =  f(u). \leqno{(\star)}
$$
\bigskip

{\bf Proposition 1.3.1.} {\it Toute section rencontre la section nulle.}\\

Il s'agit de voir qu'une application continue $f : \U \longrightarrow \C$ v\'erifiant la propri\'et\'e (*) prend la valeur~$0$. En faisant le produit membres \`a membres de la relation (*) pour $u$ et pour  son conjugu\'e $\bar{u}$, on trouve, puisque $u\bar{u} = 1$, 
$$
\overline{f(u)f(\bar{u})} = f(u)f(\bar{u}).
$$
Autrement dit, en posant $g(u) = f(u)f(\bar{u})$, on d\'efinit une application continue $g : \U \longrightarrow \R$, et il faut v\'erifier que  l'ensemble $g(\U) \subset \R$ contient l'origine. Comme $\U$ est un espace connexe et compact, $g(\U)$  est un intervalle ferm\'e de $\R$. La relation (*), avec $u = 1$ montre que $f(1)$ est r\'eel, donc que $g(1) \geq 0$ ; pour $u = -1$, on d\'eduit de cette m\^eme relation que $\overline{f(-1)} = -f(-1)$, donc que $f(-1)$  est imaginaire pur ; cela montre que $g(-1) \leq 0$. D'o\`u le r\'esultat. $\Box$
\\

{\bf Conclusion 1.3.2.}\, {\it Le ruban de M\"obius n'est pas isomorphe \`a
un cylindre.}\\

Le sens \`a donner \`a l'adjectif \emph{isomorphe} est pr\'ecis\'e dans la d\'emonstration ; elle se fait  par l'absurde en supposant qu'il existe un hom\' eomorphisme
$h:\U\times \R \; \widetilde{\longrightarrow}\; \p$  compatible avec les
projections
$$
\xymatrix{
{\U \times \R} \ar[rr]^h \ar[dr]_{{\rm pr_1}} && {\p} \ar[dl]^p\\
& \U}
$$
et induisant sur chaque fibre un isomorphisme de $\R$-espaces vectoriels \;
$\R\;
\tilde{\rightarrow}\; \p_u$. Comme, pour $u$ fix\' e, l'isomorhisme $\lambda \mapsto h(u,
\lambda)$
est suppos\'e lin\' eaire, on a $h(u, \lambda) = \lambda h(u, 1)$, donc $h(u, 1) \neq 0$.
L'application $s : \U \rightarrow \p$ d\' efinie par $ s(u) = h(u, 1)$ serait
donc une
section (continue) de $p$ partout non nulle ; on a vu que c'est impossible.
\vspace{1cm}

{\Large {\bf 2. Traduction en termes alg\' ebriques}}\\

{\bf 2.1.}\; L'espace $\p$ a \' et\' e d\' efini comme l'ensemble des points de $\U\times
\C$ qui sont fixes
sous l'involution
$$
\U\times \C\, \rightarrow \, \U\times \C, \quad (u, z)
\mapsto (u, u\bar{z}) .
$$
  Il est plus commode de
travailler avec un objet
\og param\' etr\' e\fg, c'est-\`a-dire un objet d\' efini plut\^ot comme \emph{l'image} d'une
application. Pour passer d'un point de vue \`a l'autre, on
introduit
le projecteur associ\' e, soit, ici, l'application 

$$
q(u, z) \; = \; (u, {1\over 2}(u\bar{z} + z)).
$$
Une v\' erification imm\' ediate montre que $q\circ q = q$, et que ${\rm Im}(q) = \p$.
\medskip

L'application $ z \mapsto {1\over 2}(u\bar{z} + z)$ est donc un projecteur
($\R$-lin\' eaire)  du $\R$-espace vectoriel $\C$, d\' ependant du param\`etre $u$.
Donnons-en sa matrice.

Si on \' ecrit $u = x+iy,$ avec la contrainte $x^2+y^2 = 1$, et $z = a+ib$, on
trouve
${1\over 2}(u\bar{z} + z) = {1\over 2}[((1+x)a + yb) +i(ya+(1-x)b)]$. Dans
l'espace
fibre ($= \C$) au dessus de $u = x+iy$, la matrice de
$q$ est donc
$$
\boxed{ Q \; = \; {1\over 2}{\begin{pmatrix} 1+x& y\\ y& 1-x\\ \end{pmatrix}}}
$$
\bigskip

{\bf 2.2.} \; Posons $ A = \R[X, Y]/(X^2+Y^2-1)$.
\medskip

Il faut voir $A$ comme l'anneau des fonctions continues $a : \U \rightarrow
\R$ qui
sont polynomiales au sens suivant : il existe un polyn\^ome $F(X, Y) \in
\R[X, Y]$ tel
que, pour $u = \xi +i\eta \in \U$, on ait $a(u) = F(\xi, \eta)$. Bien
\' evidemment $F$
n'est d\' efini par $a$ qu'\`a un multiple pr\`es de $X^2+Y^2-1$, et c'est
pourquoi l'on
passe au quotient.\\

D\' esignons d\' esormais par $x$ et $y$ les classes de $X$ et de $Y$ dans l'anneau
quotient $A$, et consid\' erons la matrice  \`a coefficients dans $A$
$$
Q \; = \; {1\over 2}{\begin{pmatrix} 1+x& y\\ y& 1-x\\ \end{pmatrix}}
$$
 Le th\' eor\`eme de Hamilton-Cayley
s'\' ecrit
$$
Q^2 - {\rm Tr}(Q)Q + {\rm det}(Q) = 0.
$$
Comme ${\rm Tr}(Q) = 1$ et ${\rm det}(Q) = 0$, on a
$$
Q^2 \; = \; Q.
$$
Cette matrice d\' efinit une application $A$-lin\' eaire $Q : A^2 \; \rightarrow
\; A^2$.
On pose
$$
M \; = \; {\rm Im}(Q) \; = \; {\rm Ker}(1-Q).
$$

\noindent {\bf {\it Ce module $M$ est l'analogue alg\' ebrique du ruban de M\"obius.}}
\bigskip

C'est, en fait, bien plus qu'une analogie : suivant une id\'ee  due \`a {\sc Gelfand} et maintenant classique, on \'etablit une correspondance rigoureuse entre les points de vue alg\'ebrique et g\'eom\'etrique : on reconstitue d'abord  $\U$ \`a partir de la $\R$-alg\`ebre $A$ gr\^ace \`a une bijection
$$
{\rm Hom}_{\R-Alg}(A, \R) \quad \widetilde{\longrightarrow}\quad \U.
$$
Cette application associe \`a un morphisme de $\R$-alg\`ebres $A \rightarrow \R$ le couple $(\xi, \eta) \in \R^2$ form\'e des images par ce morphisme des \'el\'ements  $x, y \in A $ ; l'\'egalit\'e $x^2+y^2 = 1$ dans $A$ se propage en la relation $\xi^2+\eta^2 = 1$ entre ces nombres r\'eels ; ils d\'efinissent donc un \'el\'ement $ u = \xi + i\eta \in \U$ (On peut aussi reconstituer la topologie de $\U$ \`a partir de  $A$). Ainsi l'ensemble, et m\^eme l'espace topologique, $\U$ est-il d\'etermin\'e par l'alg\`ebre $A$ des applications polynomiales de $\U$ dans $\R$.

Montrons ensuite comment associer \`a un \' el\' ement de $M$ une
section (polynomiale) de l'application $p : \p \rightarrow \U$. 

\`A un \' el\' ement ${a
\choose b} \in
A^2$, c'est-\`a-dire \`a un couple d'\' el\' ements de $A$, vus comme applications
polynomiales
de $\U$ dans \R, on associe l'application $a+ib: \U \rightarrow \C$; en particulier, les \'el\'ements $x$ et $y$ de $A$ doivent \^etre vus commes les fonctions coordonn\'ees, si bien que l'application $x+iy : \U \rightarrow \C$ n'est autre que l'injection $\U \subset \C$.

 Si ${a\choose b}$ est dans $M = {\rm Ker}(1-Q)$, alors $Q{a\choose b} =
{a\choose b}$, ce qui s'\' ecrit aussi, comme on le v\'erifie imm\'ediatement, 
$$
\left( {x\atop \, y}{y\atop \,-x}\right) {a\choose b} \; = \; {a\choose b}.
$$
Si on utilise l'\' ecriture complexe, cela donne 
\medskip

$$
(x+iy)(a-ib) = a+ib . \leqno{2.2.1.}
$$
\medskip

Cette relation sera r\'einterpr\'et\'ee au \S 10.7 comme un cocycle, mais il s'agit seulement ici de d\'efinir, \`a partir d'elle, une section polynomiale de $p$. Or, si on pose $f = a+ib$,  la relation 2.2.1 s'\'ecrit aussi :

$$
{\rm pour~ tout}\;   u \in \U, \; {\rm on~ a}\quad u.\overline{f(u)} = f(u).
$$
Bref, on a $(u, f(u)) \in \p$, et l'application $u \mapsto (u,
f(u))$ est la
section annonc\' ee.\\

{\bf 2.3.}\; Montrons que le $A$-module  $M$ n'est pas isomorphe \`a $A$,
contrairement \`a ce
que tout laisse penser (Au \S 4, ce module sera qualifi\'e d'inversible).

Si on pose $M' = {\rm Im}(1-Q)$, on obtient une d\' ecomposition en somme
directe de
$A$-modules
$$
M\, \oplus\, M' \; =\; A^2 .
$$
La v\' erification de cela, classique pour les espaces vectoriels, est
exactement la
m\^eme pour les modules ; elle utilise seulement le fait que l'endomorphisme
$Q$ est
idempotent. Je ne la recopie pas.\\

Soit $K$ le corps des fractions de $A$ (on verra plus bas que $A$ est
int\`egre). Notons
$Q_K$ l'endomorphisme du $K$-espace vectoriel $K^2$ de matrice $Q$, et
posons $V = {\rm
Im}(Q_K)$  et  $V' =  {\rm
Im}(1-Q_K)$. Ces $K$-espaces vectoriels donnent la d\' ecomposition
$$
V \oplus V' = K^2.
$$

Comme les matrices $Q$ et
$1-Q$ sont non
nulles, on a $\dim_K(V) = \dim_K(V') = 1$ ; $V$ et $V'$ sont, par suite,
isomorphes \`a
$K$ et ne peuvent donc pas contenir de partie libre sur $A$ ayant 2 \' el\' ements.
Enfin, l'inclusion
$A^2
\subset K^2$ entra\^ine les suivantes :
$M = Q(A^2) \subset Q(K^2) = V$, et  $M' \subset V'$.

Le module $M = {\rm Im}(Q)$ contient les \' el\' ements $Q{1\choose 0} =
{1+x\choose y}$\; et \; $Q{0\choose 1} = {y\choose 1-x}$. Si $M$ admettait
une base, elle serait r\' eduite \`a un \' el\' ement puisque $M$ est un sous-module
de $V$ ;
notons ${a\choose b}$ ce g\' en\' erateur libre suppos\'e ; il existerait des \' el\' ements
$c, d
\in A$ tels que ${1+x\choose y} = c{a\choose b}$, et ${y\choose 1-x} =
d{a\choose
b}$ ; en
particulier, on aurait  $ 1+x = ca$  et  $1-x = db$, donc $2 = ca + db$.
Mais alors
les applications $a$ et $b$ ne pourraient s'annuler simultan\' ement, et la
section $u
\mapsto (u, f(u))$ associ\' ee, comme ci-dessus, \`a ${a \choose b}$ serait
partout non nulle. On a vu que c'est impossible (1.3.1).
\medskip

Si l'anneau  $A$ \' etait principal, $M$ et
$M'$ seraient libres (comme sous-modules du module libre $A^2$), donc
libres de rang  $1$ comme sous-modules de $V$ et $V'$, et on pourrait encore
conclure que  $M$ est isomorphe \`a $A$. Cela montre g\'eom\'etriquement que $A$ n'est pas principal.
\vspace{1cm}

{\Large {\bf 3. Alg\`ebre : l'anneau}\; $A = \R[X,Y]/(X^2+Y^2-1)$}
 
\vspace{0,5cm}

Ce paragraphe aborde plusieurs aspects de la non factorialit\' e de $A$ (on explique en 3.6 pourquoi $A$ serait principal s'il \'etait factoriel). On
montre, en particulier,
que l'id\' eal
$(1+x, y)$ n'est pas principal. Ces digressions peuvent \^etre lues comme des commentaires, anticip\'es, au \S 7.\\

{\bf 3.1\; Le morphisme} \,$\R[X]\, \rightarrow \, A$.\\

Pour \' etablir certaines propri\' et\' es de l'anneau $A = \R[X,Y]/(X^2+Y^2-1)$ on
utilise  ici le morphisme d'inclusion $\R[X] \, \longrightarrow A$ qui est l'analogue alg\' ebrique de la projection du cercle unit\' e \U \, sur l'axe
des $X$.
Toute fonction polynomiale sur cet axe, i.e. tout polyn\^ome $F(X)$, fournit,
en la
composant avec la projection (la partie r\' eelle), une application
polynomiale d\' efinie
sur
$\U$, \`a savoir
$u \mapsto F(\Re u)$.\\

Notons, comme plus haut, par $x$ et $y$ les classes dans $A$ de $X$ et de
$Y$. En
\' ecrivant
$$A = \R[X][Y]/(Y^2 -(1-X^2)),$$
 on fait appara\^itre $A$ comme l'anneau obtenu
en adjoignant \`a $\R[X]$  un \'el\'ement $y$
 qui est une racine carr\' ee de $1-X^2$ ; cela montre d\' eja
que {\it tout
\' el\' ement de $A$ s'\' ecrit de fa\c{c}on unique sous la forme $p(x) + yq(x)$}.\\

Par ailleurs, le crit\`ere d'Eisenstein, appliqu\' e avec l'\' el\' ement premier $X-1$ de
$\R[X]$ entra\^ine que $Y^2 -(1-X^2)$ est irr\' eductible, donc que $A$ est
\emph{int\`egre}.\\

Il y a ici une petite difficult\'e typographique : la variable $X$ de $\R[X]$ est \'ecrite $x$ dans $A$, bien qu'elle reste trancendante sur $\R$ puisque le morphisme $\R[X] \rightarrow A$ est injectif ; il serait peut-\^etre logique - mais certainement disgracieux - d'\'ecrire $p(X)+yq(X)$.
\medskip
 
On introduit l'application norme (d\' eja utilis\' ee implicitement dans la d\' emonstration de 1.3.1.)
$$
\xymatrix{
{\R[X]\;}\ar[r] &A \ar@/_/[l]_N
}
$$
Soit $\sigma : A \rightarrow A$
l'automorphisme de
$\R[X]$-alg\`ebres d\' efini par $\sigma(y) = -y$ (C'est bien un automorphisme
puisque
$-y$ est une racine de  $Y^2 -(1-X^2)$). Si on pense \`a un \'el\'ement $a \in A$ comme \`a une application polynomiale $a : \U \rightarrow \R, \; u \mapsto a(u)$, alors $\sigma(a)$ est l'application $u\mapsto a(\bar{u})$. Posons, pour $a \in A$,
$$
N(a) = a\sigma(a).
$$
Comme $\sigma$ respecte la multiplication, il en est de m\^eme de $N$ ; on a donc
$N(ab) = N(a)N(b),$  et, si $a$ ne d\' epend pas de $y$ il est invariant par $\sigma$ et $N(a) = a^2.$

En utilisant  la base $\{1, y\}$ du $\R[X]$-module $A$, on peut \'ecrire $a = p(x)+yq(x)$, et
$$
N(p(x)+yq(x)) = p(x)^2 -y^2q(x)^2 = p(x)^2 + (x^2-1)q(x)^2.
$$
C'est un polyn\^ome en $X$, dont on peut consid\' erer le degr\' e.
\vspace{0,5cm}

{\bf Lemme 3.1.1}\; {\it Pour tout} $a\in A,\; \deg N(a) \neq 1.$\\

En effet, le polyn\^ome \`a coefficients r\' eels
$p(X)^2+(X^2-1)q(X)^2$ prend des valeurs
$\geq 0$ sur l'ouvert  $]-\infty, -1[\cup ]1, +\infty[ \, \subset \R$, ce qu'un polyn\^ome de degr\' e 1 ne
peut faire.
\vspace{0,7cm}

{\bf 3.2. Un \' el\' ement de $A$ irr\' eductible et non premier}\\

Le lemme pr\' ec\' edent
implique imm\' ediatement que
$y$ est un \' el\' ement
irr\' eductible dans $A$ : en effet, une \' egalit\' e $y = ab$, avec $a, b \in A$,
entra\^ine que
$X^2-1 = N(y) = N(a)N(b)$ ; comme $N(a)$ et $N(b)$ sont des polyn\^omes en
$X$ de degr\' e
$\neq 1$, l'un des deux est constant, donc $a$ ou $b$ est inversible.

L'irr\' eductibilit\' e de $y$ dans $A$ est li\' ee \`a la topologie de
$\R$ : cet \'el\'ement 
$y$ n'est plus irr\' eductible dans l'anneau
$\C[X, Y]/(X^2+Y^2-1)$, car on peut le d\' ecomposer en $y = {1\over
2i}(1-x+iy)(1+x+iy)$. (voir 3.5)\\

On constate ensuite que l'id\' eal $yA$ n'est pas
premier, puisque l'anneau quotient  
$$A/yA = {\mathbf R}[X, Y]/(X^2+Y^2-1,
Y) \simeq
{\mathbf R}[X]/(X^2-1)$$
 n'est pas int\`egre (la classe $\bar{x}$ de $X$ est
distincte de
$\pm 1$ et v\' erifie la relation $(\bar{x}-1)(\bar{x}+1) = 0$ ).\\

Cela montre que l'anneau $A$ n'est pas factoriel,
et a fortiori qu'il n'est pas principal.\\

 {\bf 3.3. L'id\' eal}\; $(1+x)A+yA.$\\

Notons $I$ cet id\' eal de $A$.
\medskip

3.3.1.\; On va montrer que {\it le module
$M$ introduit en 2.2 est isomorphe \`a $I$}, ce qui justifie de consid\' erer
ici cet id\' eal. Plus pr\' ecis\' ement, on va montrer que la projection 
$${\rm pr}_1 : A^2
\rightarrow A$$
induit un isomorphisme
de $M$ sur $I$. En effet, l'id\' eal $I$  est engendr\' e
par la premi\`ere ligne de la matrice $Q$, c'est-\`a-dire par les premi\`eres
coordonn\' ees de
$Q{1\choose 0}$  et $Q{0\choose 1}$. Or, ${\rm Im}(Q) = M$ ; par suite,
$$
I \; = \; {\rm pr}_1(M).
$$
Il reste \`a montrer que ${\rm Ker}({\rm pr}_1) \cap M = 0$,
c'est-\`a-dire, puisque
${\rm Ker}({\rm pr}_1) = 0\times A$, que
l'on a
$M
\cap 0\times A = 0$. Or, si
${0\choose a} \in M,$ alors $Q{0\choose a} = {0\choose a}$, ce qui \' equivaut
\`a : $ya =
0$, et $(1-x)a = a$, d'o\`u $xa = 0$ ; mais alors $a = (x^2+y^2)a = 0$.\\

On a montr\' e en 2.3, par voie g\' eom\' etrique, que le $A$-module $M$ n'est pas
libre, cela
entra\^ine donc que l'id\' eal $I$ n'est pas principal.\\

On va retrouver ce
r\' esultat presqu'alg\' ebriquement, c'est-\`a-dire en r\' eduisant les in\' evitables
arguments topologiques au seul lemme 3.1.1 ci-dessus.
\\

3.3.2.\; Montrons d'abord que l'id\'eal $I^2$ est principal, engendr\' e par $1+x$.
\medskip

L'id\' eal $I^2$ est engendr\' e par les \' el\' ements \, $(1+x)^2, (1+x)y,$  et $y^2 =
(1+x)(1-x)$ ; mettant $1+x$ en facteur, il suffit de voir que l'id\' eal engendr\' e
par $1+x$, $y$  et  $1-x$ est \' egal \`a
$A$ ; or, il contient
$(1+x)+(1-x) = 2$, qui est un \' el\' ement inversible.\\

3.3.3.\; L'id\' eal $I$ n'est pas principal.
\medskip

Raisonnons par l'absurde en supposant que $I = aA$. On a alors $(1+x)A = I^2 =
a^2A$ ; en particulier, il existe des \' el\' ements  $b$ et $c$  dans $A$ tels
que
$ a^2 = (1+x)b$ et $1+x = a^2c$. Ces \' el\' ements $b$  et $c$ sont inverses l'un
de l'autre puisque  $1+x = (1+x)bc$, et que $A$ est
int\`egre.

Prenant les normes des deux membres de l'\' egalit\' e  $ a^2 = (1+x)b$, on trouve
$$
N(a)^2 = (1+x)^2N(b).
$$
Comme $b$ est inversible, $N(b)$ est un \' el\' ement inversible de $\R[X]$,
c'est-\`a-dire
une constante non nulle ; par suite, $N(a)$ est un polyn\^ome de degr\' e 1. C'est
impossible.
\vspace{0,7cm}

{\bf 3.4. L'analogue alg\' ebrique du rev\^etement} \;$v \mapsto v^2$\\

Dans ce paragraphe et le suivant on indique comment rendre l'id\' eal $I$
principal --- et le module $M$ libre ---
en passant de $A$ \`a un sur-anneau un peu plus gros. 

Une premi\`ere m\' ethode
consiste
\`a alg\' ebriser le lemme 1.1.2 qui montre que l'application $v \mapsto v^2$
transforme le
ruban de M\"obius en un cylindre. En termes
des coordonn\' ees r\' eelles $(x, y)$ (de sorte que $v = x+iy$), cette
application s'\' ecrit
$$
(x, y) \longmapsto (x^2-y^2, 2xy).
$$
Cela indique ce qu'il faut faire.

\noindent Pour d\' efinir le morphisme de $\R$-alg\`ebres $\alpha : A\rightarrow A$
correspondant \`a
$v\mapsto v^2$, il est plus clair d'\' ecrire, dans le second anneau,
l'anneau \og but\fg ,\,
$\xi$ et
$\eta$,
\`a la place de $x$ et de $y$. On pose alors
$$
\alpha (x) = \xi^2 - \eta^2,\qquad \alpha (y) = 2\xi \eta.
$$
Cela d\' efinit bien un morphisme de $\R$-alg\`ebres puisque l'image de
$x^2+y^2-1$ est $(\xi^2-\eta^2)^2 + 4\xi^2\eta^2 - 1 = (\xi^2+\eta^2)^2 -
1 = 0$. L'id\' eal engendr\' e par
$\alpha(I)$ est engendr\' e par
$\alpha(1+x) = 1+\xi^2 -
\eta^2 =
2\xi^2$ et $\alpha(y) = 2\xi \eta$ ; comme l'id\' eal engendr\' e par $\xi$ et $\eta$ est
\' egal \`a $A$
(puisqu'il contient $1 = \xi^2+\eta^2$), on voit que
$$
\alpha(I)A = \xi A.
$$
On aurait pu aussi remarquer que la matrice $\alpha (Q)$ est
$$
{\begin{pmatrix} \xi^2& \xi\eta\\ \xi\eta& \eta^2\\ \end{pmatrix}} = {\begin{pmatrix} \xi \\ \eta\\ \end{pmatrix}}.{\begin{pmatrix} \xi & \eta  \end{pmatrix}}
$$
Par suite, ${\rm Im}(\alpha(Q))$   est visiblement le module libre engendr\'e par ${\begin{pmatrix} \xi \\ \eta\\ \end{pmatrix}} \in A^2$.
\vspace{0,5cm}

{\bf 3.5. Complexifier  $A$}\\

3.5.1.\; Consid\' erons l'anneau  $A$ comme un sous-anneau de
$B = \C[X, Y]/(X^2+Y^2-1)$.

L'anneau $B$ est principal. En effet, l'isomorphisme de \C-alg\`ebres
$$
\C[X, Y] \longrightarrow \C[Z, Z'], \qquad X\mapsto {1\over 2}(Z+Z'),\; Y
\mapsto
{1\over 2i}(Z-Z'),
$$
donne par passage aux quotients un isomorphisme
$$
\C[X, Y]/(X^2+Y^2-1) \quad \widetilde{\longrightarrow}\quad \C[Z, Z']/(ZZ'-1).
$$
Ainsi $B$
est isomorphe \`a l'anneau de fractions  $\C[Z]_Z$, lequel est principal.

D'ailleurs, il est clair que l'\' el\' ement $b = 1+x+iy \in IB$ engendre cet
id\' eal puisque
$$
1+x = {1\over 2}(1+x-iy)(1+x+iy),\qquad {\mathrm et}\qquad y = {1\over
2i}(1-x+iy)(1+x+iy).
$$
\vspace{0,5cm}

3.5.2. On va pr\' eciser les relations entre $A$ et $B$. Notons que $B =
A[i]$ (ce qui justifie le titre du paragraphe) : en effet,
chaque \' el\' ement de
$B$ s'\' ecrit de fa\c{c}on unique sous la forme $p(x) +yq(x)$, o\`u $p$ et $q$
sont des
polyn\^omes \`a coefficients complexes ; on peut donc les d\' ecomposer en $p =
p'+ip''$,
et  $q = q'+iq''$, o\`u les quatres polyn\^omes \' ecrits sont \`a coefficients r\' eels
; mais alors, $p+yq = (p'+yq') + i(p''+yq'')$.\\

On peut donc d\' efinir un automorphisme de conjugaison $b \mapsto \bar{b}$ : il 
ne porte que sur les coefficients des divers polyn\^omes, et qui laisse
invariants $x$ et $y$ ; avec cette notation, on a pour tout $b \in B$,
l'\' equivalence $ b \in A  \Leftrightarrow  b = \bar{b}$.

Introduisons les corps des fractions $K$ et $L$, respectivement de $A$ et de 
$B$ ; on a le diagramme commutatif de morphismes d'inclusion :
$$
\begin{CD}
A @>>> B\\
@VVV  @VVV\\
K @>>> L
\end{CD}
$$
Montrons l'\' egalit\' e \; $K \cap B \; = \; A$.

\noindent Il est - ou devrait \^etre - clair que $K[i] = L$, donc que $K$ est le
sous-corps de $L$ form\' e des \' el\' ements invariants sous la conjugaison invoqu\' ee
ci-dessus ; si un tel \' el\' ement est aussi dans $B$, alors il est dans~$A$.

Cela implique, en particulier, que $A$ est {\it int\' egralement clos},
autrement dit qu'un \' el\' ement de $K$ qui est racine d'un polyn\^ome unitaire
$F(T) \in A[T]$ est dans $A$ : en effet, un tel \' el\' ement, vu dans $L$ est
dans $B$ puisque $B$ est principal, donc int\' egralement clos. On n'utilisera pas cette propri\'et\'e.
\vspace{1cm}

{\bf Appendice 3.6.}\; {\bf  Sur certains anneaux factoriels}\\

 Un anneau principal est factoriel, et
il existe des anneaux factoriels qui ne sont pas principaux, par
exemple $\mathbf{Z}[X]$.

Pour clarifier le \S 3, on montre que, dans la situation rencontr\'ee,  un anneau factoriel est n\' ecessairement
principal.\\

{\bf Proposition 3.6.1.}\; {\it Soient $R$ un anneau principal contenu dans un
anneau $A$ ; on
suppose que
$A$ est un
$R$-module libre de rang fini. Alors, si $A$ est factoriel, il est
principal.}\\

Un id\' eal $I$ de $A$ est, en particulier, un sous-$R$-module du $R$-module
libre de type
fini $A$ ; c'est donc un $R$-module (libre et) de type fini ;  en
particulier $I$ est un
id\' eal de $A$ de type fini. Pour montrer qu'il est principal, on peut donc
se ramener, par
r\' ecurrence, au cas o\`u il est engendr\' e par deux \' el\' ements $a$ et $b$ ; soit
$d$ un
plus grand diviseur commun ($A$ est suppos\' e factoriel) ; \' ecrivons $a
=da'$  et
$b=db'$, de sorte que $a'$ et $b'$ n'ont pas de diviseur commun non
inversible ; on a
$aA+bA = d(a'A+b'A)$.

On est donc ramen\'e \`a d\' emontrer que si
deux \' el\' ements
$a$ et
$b$ de
$A$ n'ont pas de diviseur commun non inversible, alors la relation de B\'ezout est v\'erifi\'ee : $aA+bA = A$ ; elle \'equivaut \`a l'inversibilit\'e de $a$ modulo $b$ (c'est-\`a-dire \`a l'inversibilit\'e de l'image de $a$ dans l'anneau $A/bA$).

Si $a$ est inversible dans $A$, on a fini ; sinon, $a$ est un produit
d'\' el\' ements
irr\' eductibles de $A$. Or, si l'on
peut d\' ecomposer
$a$ en un produit
$a = a'a''$ de deux \'el\'ements $a'$  et  $a''$, inversibles modulo $b$, alors $a$ est lui-m\^eme inversible modulo $b$. Par
r\' ecurrence
sur le nombre de facteurs irr\' eductibles de $a$, il suffit donc de traiter
le cas o\`u
$a$ est irr\' eductible ; l'hypoth\`ese sur
$a$ et $b$ signifie alors que $a$ ne divise pas $b$, c'est-\`a-dire que $b
\notin aA$,
et il faut conclure que $aA +bA = A$ ; autrement dit, il faut montrer que
l'id\' eal $aA$
est \emph{maximal}. En changeant la notation, et en utilisant le fait que
dans un
anneau factoriel, un \' el\' ement irr\' eductible est premier, on est ramen\' e \`a
d\' emontrer ceci :
 {\it sous les hypoth\`eses de la proposition, si
$p$ engendre un id\'eal premier non nul de
$A$, alors cet id\' eal
$pA$ est maximal, i.e $A/pA$ est un corps.}

On remarquera que la factorialit\'e de $A$ n'est plus utilis\'ee dans la suite ; par contre le fait que $A$ soit fini et libre sur le sous-anneau principal $R$ est essentiel.  

Montrons d'abord que l'id\' eal $pA \cap R$ est non nul. L'application $u :
A\rightarrow A,
\; a
\mapsto pa$ est un endomorphisme $R$-lin\' eaire du $R$-module libre $A$. Le th\'eor\`eme de Hamilton-Cayley montre qu'il existe un polyn\^ome unitaire $F(T) \in R[T]$ tel que $F(p) = 0$. Choisissons un polyn\^ome unitaire $F \in R[T]$, annulant $p$ et de degr\'e minimum $d$, soit  
$$
F(T)\; =\; X^d +r_{d-1}X^{d-1} + \cdots +r_{0}
$$
Le terme constant $r_{0}$  est non nul, sinon $F(T) = TF_{1}(T)$, et l'int\'egrit\'e de $A$ montre que $F_{1}(p) = 0$, contrairement \`a la minimalit\'e du degr\'e. Il est clair que $0 \neq r_{0} \in R \cap pA$.

En passant aux quotients, on obtient un homomorphisme \emph{injectif} d'anneaux
$$
R/pA\cap R\; \longrightarrow \; A/pA .
$$
Comme $pA$ est premier, l'anneau $A/pA$ est int\`egre, donc son
sous-anneau
$R/pA\cap R$ l'est aussi ; l'id\' eal $pA\cap R$ est donc premier, et par
suite maximal puisqu'il est
non nul et que $R$ est principal ; bref, $R/pA\cap R$ est un corps, et
l'anneau int\`egre
$A/pA$ appara\^it comme un espace vectoriel de dimension finie sur ce corps ;
on en
d\' eduit que $A/pA$ est un corps, en invoquant le lemme classique suivant.\\

{\bf Lemme 3.6.2.}\; {\it Un anneau int\`egre $S$ qui est un espace vectoriel de dimension finie sur un
sous-corps $K$, est lui-m\^eme un corps.}\\

En effet, si $s\in S$ est un \'el\'ement non nul, la multiplication par $s$ est un endomorphisme
injectif de $S$, puisque $S$ est int\`egre, donc bijectif car $S$ est un vectoriel de dimension finie
sur le sous-corps $K$ ; l'\' el\' ement
unit\' e $1$ est donc dans l'image de cet endomorphisme ; ainsi $s$ admet un
inverse.
\vspace{1cm}

{\Large {\bf 4. Modules inversibles}}
\bigskip

Il y a, au moins, trois d\'efinitions possibles pour la notion de module inversible. On donne d'abord ici la plus concr\`ete et qui demande le moins de pr\'eliminaires. Les deux autres seront introduites dans les \S\S 5 et 8, ainsi que leur \'equivalence.
\bigskip

{\bf  4.1. \; Un module inversible est un module projectif de rang 1}
\bigskip

4.1.1\; Un $A$-module $L$ est dit \emph{projectif} (de type fini) s'il est facteur direct d'un $A$-module libre de type fini, autrement dit, s'il existe un $A$-module $L'$ et un isomorphisme $ L\oplus L' \simeq A^n$. Il revient au m\^eme de dire qu'il existe une application $A$-lin\'eaire $f : A^n \rightarrow A^n$ qui est un projecteur ($f^2 = f$), et un isomorphisme
$$
{\rm Im}(f) \; \simeq \; L.
$$
L'endomorphisme $f$ se factorise en $ f = u\circ v$, o\`u $u : L \rightarrow A^n$ est l'injection canonique, et $v : A^n \rightarrow L$ est la surjection d\'eduite de $f$. Comme $f^2 = f$, on a $uvuv = uv$, mais $v$ est surjectif et $u$ est injectif ; on a donc  deux applications
$$
L\; \stackrel{u}{\longrightarrow}\; A^n\; \stackrel{v}{\longrightarrow}\; L \quad {\rm telles ~que }\quad vu = {\rm Id}_{L}.
$$
\medskip

4.1.2 \; Pour un espace vectoriel (de dimension finie), le mot  \emph{rang} est synonyme de \emph{dimension}, et on emploie indiff\'eremment l'un ou l'autre. 

Soit $\pp$ un id\'eal premier de $A$, de sorte que l'anneau quotient $A/\pp$ est int\`egre, et admet donc un corps des fractions que l'on note $\kappa(\pp)$.

Soit $L$ un $A$-module de type fini, et soit $\pp$ un id\'eal premier ; le quotient $L/\pp L$  est un module sur l'anneau int\`egre $A/\pp$, et  son localis\'e $(L/\pp L)_{\pp}$ est un $\kappa(\pp)$-espace vectoriel de dimension finie ; on pose
$$
{\rm rang}_{\pp}(L) \; = \; {\rm dim}_{\kappa(\pp)}((L/\pp L)_{\pp}).
$$
Un $A$-module de type fini $L$  est dit \emph{de rang $n$} si l'application $\, \pp \mapsto  {\rm rang}_{\pp}(L) $  est constante de valeur ~$n$.\medskip

\noindent Le rang d'un module projectif de type fini se trouve \^etre \'egal au rang d'un projecteur associ\'e ; plus pr\'ecisement  :\\

{\bf Lemme  4.1.3.}\;  {\it Soit $f = f^2$ un projecteur de $A^n$, et soit $L = {\rm Im}(f)$ son image. Soit $\pp$ un id\'eal premier de $A$ ; notons $k = \kappa(\pp) $ le corps des fractions de $A/\pp$, et affectons d'une barre les $k$-endomorphismes obtenus par r\'eduction modulo $\pp$ et passage au corps des fractions. Alors, le  $k$-espace vectoriel $(L/\pp L)_{\pp}$ est canoniquement isomorphe \`a l'image de l'endomorphisme $\bar{f} : k^n \rightarrow k^n$ ; en particulier, la dimension  de $(L/\pp L)_{\pp}$, soit ${\rm rang}_{\pp}(L)$,  est \'egale au rang de l'endomorphisme $\bar{f}$.}
\medskip

Consid\'erons le diagramme commutatif suivant
$$
\begin{CD}
A^n @>{v}>> L @>{u}>> A^n @>{v}>> L\\
@VVV  @VVV  @VVV @ VVV\\
k^n @>>{\bar{v}}> (L/\pp L)_{\pp} @>>{\bar{u}}> k^n @>>{\bar{v}}> (L/\pp L)_{\pp}
\end{CD}
$$
On va voir que $\bar{u}$ \'etablit un isomorphisme de $(L/\pp L)_{\pp}$ sur ${\rm Im}(\bar{f})$.
Les deux carr\'es de gauche montrent ceci : comme $v$ est surjectif, $\bar{v}$ l'est aussi ; donc ${\rm Im}(\bar{f}) = {\rm Im}(\bar{u}\bar{v}) =  {\rm Im}(\bar{u})$ ; il reste donc \`a v\'erifier que $\bar{u}$ est injectif. Mais, en regardant les deux carr\'es de droite, on voit que  $\bar{v}\bar{u} = {\rm Id}$, donc que $\bar{u}$ est injectif.$\Box$\\

{\bf D\'efinition 4.1.4}\; {\it Un $A$-module est dit \emph{inversible} s'il est projectif, de type fini et de rang $1$.}
\medskip

L'anneau $A$ est lui-m\^eme un $A$-module inversible ; un module inversible et libre est isomorphe \`a $A$; on dit alors  volontiers qu'un tel module inversible est \emph{trivial}. En g\'en\'eral, montrer qu'un module inversible \emph{n'est pas} libre est chose malais\'ee.
\bigskip

{\bf Proposition 4.1.5.}\; {\it Pour un endomorphisme $f$ de $A^2$, la condition : \^etre un projecteur dont l'image est de rang 1, \'equivaut \`a}
$$
{\rm Tr}(f) = 1,\quad et \quad  {\rm det}(f) = 0. \leqno{(\star)}
$$
\medskip

La d\'emonstration de la proposition utilise le th\'eor\`eme de Hamilton-Cayley, qui s'\'ecrit ici
$$
f^2 - {\rm Tr}(f)f + {\rm det}(f) = 0.
$$
Si les conditions $(\star)$ sont v\'erifi\'ees, alors $f^2 = f$. Posons $L = {\rm Im}(f)$ ; il s'agit de montrer que $L$ est de rang 1. Consid\'erons donc un id\'eal premier $\pp$ de $A$ ; d'apr\`es le lemme 4.1.3, dont nous gardons les notations,  
il faut voir que le rang du projecteur $\bar{f}$ du $k$-espace vectoriel $k^2$,  est \'egal \`a 1. Or, $\bar{f}$ est non nul - donc de rang $\geq 1$ -  puisque sa trace est 1, et il n'est pas un isomorphisme puisque son d\'eterminant est nul.  \medskip

Avant de d\'emontrer la r\'eciproque, rappelons le 
\medskip

{\bf Lemme 4.1.6}\; {\it Soit $e =e^2$ un idempotent d'un anneau (commutatif) $A$. Si $e$ n'est contenu dans aucun id\'eal premier   de $A$ (resp. s'il est contenu dans tous), alors $e = 1$ (resp. $e = 0$).}
\medskip

Si  $e$ n'est contenu dans aucun id\'eal premier, il est inversible ; mais $e(1-e) = 0$, donc $1-e =0$ ; l'assertion parall\`ele se d\'eduit de la premi\`ere en rempla\c{c}ant $e$ par $1-e$. $\Box$ \medskip

Pour achever la d\'emonstration de la proposition, montrons qu'un projecteur $f$ de $A^2$ dont l'image est un module de rang 1 v\'erifie les relations $(\star)$. D'apr\`es le lemme 4.1.3, pour tout id\'eal premier  $\pp$, le rang de l'endomorphisme $\bar{f}$ est \'egal \`a $1$ ; en particulier, son d\'eterminant est nul ; bref, ${\rm det}(f)$ est contenu dans tous les id\'eaux premiers de $A$ ; mais, par ailleurs, c'est un idempotent, tout comme $f$ ; donc ${\rm det}(f) = 0$, d'apr\`es  4.1.6.

Le th\'eor\`eme de Hamilton-Cayley, et la relation $f^2 = f$ entra\^inent alors que 
$$
(1-{\rm Tr}(f)) f  = 0. \leqno{(\star \star)}
$$
En prenant la trace, qui est $A$-lin\'eaire, on en d\'eduit que ${\rm Tr}(f)$ est un idempotent de $A$ ; s'il \'etait contenu dans un id\'eal premier $\pp$, la relation $(\star \star)$ montre qu'on aurait $\bar{f} = 0$, ce qui est impossible d'apr\`es le lemme 4.1.3.  puisque $L$ est suppos\'e de rang 1. Ainsi, l'idempotent ${\rm Tr}(f)$ n'est contenu dans aucun id\'eal premier ; il est donc \'egal \`a 1 d'apr\`es 4.1.7, et la proposition est d\'emontr\'ee. $\Box$
\bigskip

Le r\'esultat suivant sera utilis\'e plus bas.
\medskip

{\bf Lemme 4.1.7.} {\it Soit $L$ un module inversible et $\alpha : L \rightarrow A$ une application lin\'eaire surjective. Alors $\alpha$ est un isomorphisme}
\medskip

Tout $x \in L$ tel que $\alpha(x) = 1$ engendre un suppl\'ementaire dans $L$ de $M = {\rm Ker}(\alpha)$ ; pour tout id\'eal premier $\pp$ de $A$, l'application entre espaces vectoriels de rang 1, $\bar{\alpha} : (L/\pp L)_{\pp} \rightarrow (A/\pp)_{\pp} = \kappa(\pp)$ est surjective, donc bijective ; par suite, pour tout id\'eal premier $\pp$, on a $(M/\pp M)_{\pp} = 0$, et il faut en conclure que $M = 0$. Or, comme $M$ est facteur direct de $L$, c'est aussi un module projectif ; il existe donc un projecteur $g : A^m \rightarrow A^m$ d'image isomorphe \`a $M$. Le lemme 4.1.3  montre que pour tout $\pp$, le rang de $\bar{g} : \kappa(\pp)^m  \rightarrow  \kappa(\pp)^m$ est nul, donc que $1-\bar{g} = 1$ ;  on d\'eduit de 4.1.6, que l'on a ${\rm det}(1-g) = 1$, donc que $1-g$ est un isomorphisme ; mais $g(1-g) = 0$ ; donc  $g = 0$.$\Box$ 
\bigskip

 On aura remarqu\'e que la matrice $Q$ de 2.2 d\'efinit un projecteur de rang 1, dont l'image est le module $M$ que l'on montre en 2.3 \^etre inversible et non libre. C'est un cas particulier de la construction qui suit.
 \bigskip

\noindent {\bf 4.2\; Un exemple tr\`es g\'en\'eral de module inversible  et non libre}
\bigskip
 
 On pose $A = \Z [X, Y, Z]/(X^2-X+YZ)$, et on d\'esigne par $x, y$ et $z$ les classes de $X, Y$ et $Z$. On consid\`ere l'endomorphisme de $A^2$ d\'efini par $$
 f = {\begin{pmatrix} x& z\\ y& 1-x\\ \end{pmatrix}} .
 $$
D'apr\`es 4.1.5,  $L = {\rm Im}(f)$ est un $A$-module inversible.
\medskip

\noindent \emph{Montrons qu'il n'est pas libre.}
\medskip

Notons $R = \Z[y, z]$ le sous-anneau de $A$ engendr\'e par $y$ et $z$ ; ces \'el\'ements sont alg\'ebriquement ind\'ependants sur $\Z$ ; par suite $R$ est un anneau factoriel, ce que n'est pas $A$. On se ram\`ene dans $R$ en utilisant une norme que l'on d\'efinit comme suit : le $R$-module $A$ est libre de base $\{1, x\}$ ; tenant compte de la relation $ x^2 = x - yz$, on voit que la mutiplication dans $A$ par l'\'el\'ement $\alpha + \beta x$, avec $\alpha, \beta \in R$,  a pour matrice, sur la base $\{1, x\}$,
$$
 {\begin{pmatrix}\alpha & -\beta yz\\ \beta &\alpha + \beta\\ \end{pmatrix}}  .
 $$
On d\'efinit  l'application norme comme le d\'eterminant de cette matrice, soit
$${\sf N}
: A \longrightarrow R,\quad \alpha+\beta x \longmapsto \alpha^2 + \alpha \beta + \beta^2 yz.
$$
On a, en particulier, ${\sf N}(x) = {\sf N}(1-x) = yz $. La norme est une application multiplicative. 
\medskip

\noindent On utilisera la remarque suivante :
\medskip

{\it Une \'egalit\'e de la forme
$
{\sf N}(\alpha + \beta x) = ny,
$
est impossible si  $n$ est un entier non nul.} 
\medskip

Elle impliquerait, en effet, l'\'egalit\'e suivante entre polyn\^omes de $\Z[y, z]$ : 
$$
4{\sf N}(\alpha + \beta x) = (2\alpha + \beta)^2 + \beta^2(4yz -1) = 4ny .
$$ 
On en d\'eduirait que pour tout couple de r\'eels $y, z$ tels que $4yz - 1 \geq 0$, on devrait avoir $4ny \geq 0$ ; or, si $(y, z)$ est un tel couple, on a aussi $4(-y)(-z) -1 \geq 0$, ce qui conduit \`a une contradiction. $\Box$
\medskip

Supposons  que le module $L = {\rm Im}(f)$ soit libre ; cela permet de l'identifier \`a $A$, et d'\'ecrire les applications associ\'ees \`a $f$, comme en 4.1.1, \,  $ v : A^2 \rightarrow L$, et $u : L \rightarrow A^2$, sous forme de matrices \`a coefficients dans $A$, soit $v = ( v_{1}\, v_{2})$,\;et \; $u = {u_{1}\choose u_{2}}$ ; la relation $ f = uv$ s'\'ecrit alors
$$
{\begin{pmatrix} x& z\\ y& 1-x\\ \end{pmatrix}}  \; = \; \begin{pmatrix}u_{1}\\ u_{2}\\ \end{pmatrix}. \begin{pmatrix} v_{1}&v_{2}\\ \end{pmatrix} .
$$
On a donc les \'egalit\'es suivantes entre \'el\'ements de $A$
$$
\begin{array}{c}
x\, =\, u_{1}v_{1}\\
y\, =\, u_{2}v_{1}\\
z \, =\, u_{1}v_{2}\\
1-x\, =\, u_{2}v_{2}
\end{array}
$$
Prenant les normes, on trouve les \'egalit\'es suivantes entre polyn\^omes de $R = \Z[y, z]$ :
$$
\begin{array}{c}
yz\, =\, {\sf N}(u_{1}){\sf N}(v_{1})\\
y^2\, =\, {\sf N}(u_{2}){\sf N}(v_{1})\\
z^2 \, =\, {\sf N}(u_{1}){\sf N}(v_{2})\\
yz\, =\, {\sf N}(u_{2}){\sf N}(v_{2}).
\end{array}
$$
\'Etant factoriel, l'anneau $R$ poss\`ede des pgcd, et pgcd$(y, z) = 1$. 
La deuxi\`eme et la quatri\`eme \'egalit\'e donnent
$$
y = {\rm pgcd}(y^2, yz) = {\sf N}(u_{2}){\rm pgcd}({\sf N}(v_{1}), {\sf N}(v_{2})).
$$
La premi\`ere et la troisi\`eme
$$
z = {\rm pgcd}(yz, z^2) = {\sf N}(u_{1}){\rm pgcd}({\sf N}(v_{1}), {\sf N}(v_{2})).
$$
Ces deux relations montrent que l'\'el\'ement $\epsilon = {\rm pgcd}({\sf N}(u_{1}), {\sf N}(u_{2}))$ doit diviser $y$ et $z$ ; il est donc inversible dans $R = \Z[y, z]$, soit $\epsilon = \pm 1$. En utilisant maintenant les deux premi\`eres \'egalit\'es, on obtient 
$$
y = {\rm pgcd}(yz, y^2) = {\sf N}(v_{1}){\rm pgcd}({\sf N}(u_{1}), {\sf N}(u_{2})) = \epsilon {\sf N}(v_{1}).
$$
 Mais on a montr\'e au d\'ebut qu'une relation de la forme ${\sf N}(v_{1}) = \epsilon y$ \'etait impossible. Cette contradiction montre que $L$ n'est pas libre.
\bigskip

{\bf Exercice 4.2.1.}\; Montrer que la matrice de $1-f$ est semblable \`a la transpos\'ee de celle de $f$. En d\'eduire que dans la d\'ecomposition $L \oplus L' = A^2$ associ\'ee \`a $f$, le facteur $L'$ est isomorphe au dual de $L$.

Montrer que si $y = z$ (dans un quotient de $A$ !), alors $L$ et $L'$ sont isomorphes.
\bigskip

{\bf Explication 4.2.2.}\; ${}^{\ast}$La qualification de \emph{tr\`es g\'en\'eral} attribu\'ee \`a cet exemple peut \^etre  pr\'ecis\'ee de la fa\c{c}on suivante : pour tout anneau $B$ et tout $B$-module inversible $M,$ \emph{engendr\'e par deux \'el\'ements}, il existe un morphisme $A \rightarrow B$, et un isomorphisme de $B$-modules $B\otimes_{A}L\;  \widetilde{\longrightarrow}\; M$.${}_{\ast}$
\vspace{1cm}

{\Large {\bf 5. Modules localement isomorphes \`a $A$ }}
\bigskip

On utilise librement dans ce paragraphe la notion d'anneau de fractions ; on note $A_{s}$ l'anneau des fractions dont le d\'enominateur est une puissance de $s$.
\bigskip

 {\bf 5.1.  Modules localement isomorphes \`a $A$ }
\bigskip

L'adverbe \emph{localement}  renvoie ici \`a l'ensemble $\Sp (A)$ des id\'eaux premiers de $A$, muni de la \emph{topologie de Zariski} : rappelons  que les ouverts pour cette topologie  sont  les r\'eunions d'ensembles de la forme ${\rm D}(s)$, o\`u, pour  $s \in A$, on note ${\rm D}(s) = \Sp (A_{s})$ l'ensemble des id\'eaux premiers qui ne contiennent pas~$s$. 

Une famille $({\rm D}(s_{i}))_{i\in I}$ forme un recouvrement de $\Sp (A)$ si aucun id\'eal premier ne contient tous les $s_{i}$, c'est-\`a-dire si l'id\'eal $\sum As_{i}$  est \'egal \`a $A$ ; on dit alors que ces \'el\'ements sont \emph{\'etrangers}  (dans leur ensemble ; ne pas confondre cette notion avec celle, plus forte, de famille d'\'el\'ements deux \`a deux \'etrangers). Mais, pour un id\'eal, \^etre \'egal \`a $A$ \'equivaut \`a contenir $1$, et l'inclusion $1 \in \sum As_{i}$ ne fait intervenir qu'un nombre fini de $s_{i}$ ; par suite une famille d'\'el\'ements \'etrangers contient une sous-famille finie qui est encore form\'ee d'\'el\'ements \'etrangers ; autrement dit, $\Sp (A)$ est un espace quasi-compact.
\medskip

Sans pouvoir d\'evelopper les d\'etails qui justifieraient une telle interpr\'etation, disons que, pour un $A$-module $M$, un \'el\'ement du module de fractions  $M_{s}$ peut \^etre vu comme \emph{une section\footnote{Le mot \og section \fg peut sembler \^etre employ\'e ici dans un sens diff\'erent que dans 1.3 ; il n'en est rien. Comme il est expliqu\'e, par exemple dans {\sc EGA I} 1.3, \`a tout $A$-module est associ\'e un faisceau $\widetilde{M}$ sur l'espace ${\rm Spec}(A)$ de telle sorte qu'on ait une bijection
$$
M_{s} \simeq \Gamma({\rm D}(s), \widetilde{M}).
$$
Par ailleurs la notion de faisceau remplace celle, \'equivalente, d'espace \'etal\'e, ici au dessus de ${\rm Spec}(A)$, pour laquelle la notion de section est celle de $1.3.$. } de $M$ au dessus de l'ouvert ${\rm D}(s)$} ; c'est pourquoi on nomme souvent $M_{s}$ un \emph{localis\'e} de $M$ ; de m\^eme, le lemme suivant peut (devrait !) \^etre lu ainsi : une section localement nulle est nulle.
\bigskip

{\bf Lemme 5.1.1.}\; {\it Soit $(s_{i})_{i}$ une famille finie d'\'el\'ements \'etrangers de $A$. Pour tout $A$-module $M$, l'application $M \, \rightarrow \, \prod_{i} M_{s_{i}}$ est injective.}
\medskip

C'est une cons\'equence du fait que l'anneau $B = \prod A_{s_{i}}$  est fid\`element plat sur $A$ (cf \S 9.2), mais on peut aussi le voir de fa\c{c}on \'el\'ementaire : l'annulateur ${\rm Ann}(x)$ d'un \'el\'ement  $x$ du noyau contient donc une puissance de chaque $s_{i}$ ; comme ces \'el\'ements sont \'etrangers, l'id\'eal ${\rm Ann}(x)$ ne peut \^etre contenu dans aucun id\'eal maximal de $A$ ; on a donc ${\rm Ann}(x) = A$ ; ainsi, $1$ annule $x$ \dots .
\medskip

Ceci rappel\'e, le titre du paragraphe se pr\'ecise en l'\'enonc\'e suivant.
\bigskip

{\bf  Proposition 5.1.2.}\; {\it  Pour qu'un $A$-module $L$   soit inversible, il faut et il suffit  qu'il existe une famille (finie) $(s_{i})$ d'\'el\'ements \'etrangers et pour chaque $i$ un isomorphisme de $A_{s_{i}}$-modules $A_{s_{i}} \simeq L_{s_{i}}$.}
\bigskip

Pour montrer que la condition est n\'ecessaire, il faut trouver, pour chaque id\'eal maximal $\mathfrak{m}$, un \'el\'ement  $ t \in A - \mathfrak{m}$ et un isomorphisme $L_{t} \rightarrow A_{t}$ ; il suffit m\^eme, en vertu de 4.1.7, de trouver une telle application qui soit surjective. Comme $L$ est projectif, il existe deux applications $A$-lin\'eaires $L \stackrel{u}{\longrightarrow} A^n \stackrel{v}{\longrightarrow} L$  telles que $v\circ u = {\rm Id}_{L}$. Comme, par hypoth\`ese, $L/\mathfrak{m} L$ est de rang 1, il existe un \'el\'ement $x \in L$ tel que $ x \notin \mathfrak{m} L$. Les coordonn\'ees $(a_{1}, \ldots , a_{n})$ de $u(x)$ ne sont pas toutes dans $\mathfrak{m}$, sinon $x = v(u(x))$ serait dans $\mathfrak{m} L$ ; notons pour simplifier $t = a_{i}$ une coordonn\'ee non dans $\mathfrak{m}$, et soit ${\rm pr}_{i} : A^n \rightarrow A$ la projection sur le facteur  d'indice $i$ ; l'image de $x$ par l'application compos\'ee  $L \stackrel{u}{\longrightarrow} A^n \stackrel{{\rm pr}_{i}}{\longrightarrow} A$ est \'egale \`a $t$ ; l'application $L_{t} \rightarrow A_{t}$, obtenue par passage aux anneaux de fractions est donc surjective.

Le fait que cette condition soit suffisante peut \^etre d\'emontr\'e directement ;  mais  c'est un cas particulier du th\'eor\`eme de descente (cf  9.3), puisque, comme d\'ej\`a dit, l'anneau $B = \prod A_{s_{i}} $ est fid\`element plat sur $A$, et que l'hypoth\`ese sur $L$ signifie exactement qu'il existe un isomorphisme de $B$-modules $B\otimes_{A}L \; \widetilde{\longrightarrow}\; B$. La d\'emonstration directe suivrait pas \`a pas celle du th\'eor\`eme de descente, mais en beaucoup moins lisible parce qu'elle doit mettre en jeu des familles de multi-indices. Nous renvoyons donc au \S 9.$\Box$
\medskip

{\bf Corollaire 5.1.3.}\; {\it Soient $L$ un $A$-module inversible, et $u : L \rightarrow A$ une application $A$-lin\'eaire. Alors, pour tous $x, y \in L$, on a}
$$
u(x)y \; =\; u(y)x.
$$
\medskip

Il s'agit de v\'erifier une \'egalit\'e entre \'el\'ements de $L$ ; compte-tenu de 5.1.1 et 5.1.2, il suffit de le faire dans les localis\'es $L_{s}$ qui sont isomorphes \`a $A_{s}$ ; mais si $L$ est libre de rang 1, de base $z$, il existe des scalaires $a, b \in A$, tels que $x = az$  et $y =bz$ ; l'\'egalit\'e r\'esulte alors de la commutativit\'e de $A$. $\Box$
\medskip

{\bf Exemple 5.1.4.}\; Reprenons l'exemple, donn\'e en 4.2, du module $L$ image du projecteur de $A^2$ de matrice
$$
 {\begin{pmatrix} x& z\\ y& 1-x\\ \end{pmatrix}}
 $$
 On va construire des isomorphismes $A_{x} \simeq L_{x}$,\, et $A_{1-x} \simeq L_{1-x}$, en utilisant la relation
 $$
 x.(1-x) \; =\; yz.
 $$
 Le module $L \subset A^2$ est engendr\'e par les \'el\'ements ${x \choose y}$  et ${z \choose 1-x}$.
 Dans l'anneau $A_{x}$, l'\'el\'ement $x$ est inversible ; on peut donc \'ecrire
 $$
 {\begin{pmatrix} x\\ y \end{pmatrix}} \, = \, x{\begin{pmatrix} 1\\ y/x\\ \end{pmatrix}},\qquad {\begin{pmatrix}  z\\  1-x\\ \end{pmatrix}} \, = \, z{\begin{pmatrix} 1\\ y/x\\ \end{pmatrix}}.
 $$
 Par suite le module $L_{x} \subset A_{x}^2$ est engendr\'e par l'\'el\'ement ${1 \choose y/x}$ qui est visiblement libre. De m\^eme, dans $A_{1-x}$, on a
 $$
 {\begin{pmatrix} x\\ y\\ \end{pmatrix}}\, =\, y{\begin{pmatrix}z/(1-x)\\ 1\\ \end{pmatrix}}, \qquad {\begin{pmatrix}z\\ 1-x\\ \end{pmatrix}}\, =\, (1-x){\begin{pmatrix} z/(1-x)\\ 1\\ \end{pmatrix}}.
 $$
Par suite le module $L_{1-x} \subset A_{1-x}^2$  est engendr\'e par ${z/(1-x) \choose 1}$.
\vspace{10mm}

{\bf  5.2. Construction par recollement}
\bigskip

Commen\c{c}ons par le cas o\`u $A$ est int\`egre, ce qui simplifie beaucoup les choses.
\bigskip

{\bf Exercice 5.2.1.}\; Soient $A$ un anneau int\`egre de corps des fractions $K$. Soient $s, t \in A$ deux \'el\'ements \'etrangers et non nuls ; on identifie $A_{s}$  et $A_{t}$ \`a des sous-anneaux de $K$.

{\it i)}\; Montrer que l'on a $A_{s} \cap A_{t} = A$.

{\it ii)}\;  Soit $M$ un sous-$A$-module (quelconque) de $K$. Montrer que $M_{s} \cap M_{t} = M$.

{\it iii)}\; On consid\`ere ici l'anneau $A$ et le $A$-module inversible $L$ introduits au \S4.2. L'application compos\'ee 
$$L \; \stackrel{v}{\longrightarrow}\; A^2\; \stackrel{{\rm pr}_{1}}{\longrightarrow}\;  A$$
 est injective (cf 7.1.2), et induit un isomorphisme de $L$ sur le sous-module $M = Ax+Az \subset K$. Montrer que $M_{x} = A_{x}$, et $M_{1-x} = \frac{x}{y}.A_{1-x}$. En d\'eduire que $L$ est isomorphe au sous-module $\frac{y}{x}.A_{x} \cap A_{1-x}  \subset K$.
 
 {\it iv)}\; Revenons au cas g\'en\'eral. Soit $\omega$ un \'element inversible de l'anneau $A_{st}$. Alors
 $L =  \omega.A_{s} \cap A_{t}$ est un $A$ module inversible.
\bigskip

Ce qui suit g\'en\'eralise la question  {\it iv)}, lorsqu'on ne suppose plus que $A$ est int\`egre.
\bigskip

{\bf Lemme 5.2.2}\; {\it  Soit $A$ un anneau et soient $s, t$  deux \'el\'ements \'etrangers, de sorte que $sA + tA = A$. Consid\'erons les morphismes canoniques
$$
A_{s}\; \stackrel{\alpha}{\longrightarrow}\; A_{st}\; \stackrel{\beta}{\longleftarrow}\; A_{t}.
$$

\noindent Alors, la suite 
$$
\begin{CD}
0 @>>>  A@>>> A_{s}\times A_{t} @>{(x, y) \mapsto \alpha(x)-\beta(y)}>>  A_{st}
\end{CD}
$$  est exacte.}
\medskip

L\`a encore, on peut invoquer le fait que le morphisme $A \rightarrow A_{s}\times A_{t}$ est fid\`element plat (cf \S 9.2), mais on peut aussi v\'erifier directement l'exactitude : consid\'erons des \'el\'ements $a/s^n \in A_{s}$  et $b/t^m \in A_{t}$ dont les images dans $A_{st}$ sont \'egales ; cela signifie qu'il existe un entier $p$ tel que l'on ait, dans $A$,
$$
at^m.s^pt^p = bs^n.s^pt^p.
$$
Quitte \`a r\'ecrire $a/s^n$ sous la forme $as^p/s^{n+p}$, et de m\^ eme pour $b/t^m$, on peut supposer que l'on a
$$
at^m \; =\; bs^n .
$$
Comme $s$ et $t$ sont \'etrangers, $s^n$ et $t^m$ le sont aussi, et il existe dans $A$ des \'el\'ements $u$ et $v$ tels que $s^nu + t^mv = 1$. Posons  $c = bv+au$. On trouve, dans $A_{s}$,
$$
a/s^n \; =\; (at^mv+as^nu)/s^n \; =\; (bv+au)/1 \; = \; c/1 .
$$
De m\^eme, dans $A_{t}$, on a $b/t^m = c/1$. D'o\`u le r\'esultat.
\bigskip

{\bf Proposition 5.2.3.}\; {\it Soient $s$  et $t$  des \'el\'ements \'etrangers dans $A$, et $\omega$ un \'el\'ement inversible de l'anneau $A_{st}$. Alors le $A$-module}
$$
L \; =\; \{(\xi, \eta) \in A_{s}\times A_{t},\; \omega \alpha(\xi) = \beta(\eta)\}
$$
{\it est inversible. Il est libre si et seulement si il existe des \'el\'ements \emph{inversibles} $u \in (A_{s})^{\times}$  et $v \in (A_{t})^{\times}$ tels que $(u, v) \in L$, c'est-\`a-dire tels que $\omega \alpha(u) = \beta(v)$.}
\medskip

Par d\'efinition de $L$, on a la suite exacte
$$
\begin{CD}
0 @>>>  L@>>> A_{s}\times A_{t} @>{(\xi, \eta) \mapsto \omega\alpha(\xi)-\beta(\eta)}>>  A_{st} .
\end{CD}
$$ 
Elle reste exacte par localisation, et s'\'ecrit
$$
\begin{CD}
0 @>>>  L_{s}@>>> A_{s}\times A_{st} @>{(\xi, \eta) \mapsto \omega\alpha(\xi)- \eta}>>  A_{st}
\end{CD}
$$ 
Il est alors clair que l'\'el\'ement $(1, \omega) \in A_{s}\times A_{st}$ est dans $L_{s}$, et en constitue une base comme $A_{s}$-module. De m\^eme, l'\'el\'ement $(\omega^{-1}, 1) \in A_{st}\times A_{t}$ est une base de $L_{t}$. Ainsi, d'apr\`es 5.1.2,  $L$ est inversible.

Supposons que $L$ soit libre, engendr\'e par $(u, v) \, \in \,  A_{s}\times A_{t}$, et montrons que $u$ est inversible dans $A_{s}$. Par localisation en $s$, l'\'el\'ement $(u, \beta(v)) \in A_{s}\times A_{st}$ est une base du $A_{s}$-module $L_{s}$ ; il existe donc un \'el\'ement $a \in A_{s}$ tel que
$$
(1, \omega) = a(u, \beta(v)).
$$
Par suite, on a, dans l'anneau $A_{s}$, \, $1 = au$, et $u$ est bien inversible ; on montrerait de m\^eme que $v$ est inversible dans $A_{t}$. 

R\'eciproquement, montrons qu'un \'el\'ement  $(u, v) \, \in \,  (A_{s})^{\times} \times (A_{t})^{\times}$ tel que $\omega\alpha(u) = \beta(v)$ est une base de $L$. Or, pour tout $(\xi, \eta) \in L$, l'\'el\'ement $(\xi.u^{-1}, \eta.v^{-1}) \, \in \, A_{s}\times A_{t}$ v\'erifie la relation $\alpha(\xi.u^{-1}) = \beta(\eta.v^{-1})$ ; il provient donc de $A$, comme il est rappel\'e au d\'ebut (5.2.2), ce qui veut dire qu'il existe $a \in A$ tel que $\xi = au$  et $\eta = av$. $\Box$
\bigskip

On peut, bien entendu, construire un module inversible par recollement \`a partir d'un nombre fini quelconque d'\'el\'ements \'etrangers $s_{1}, \ldots, s_{n}$, et d'\'el\'ements inversibles $\omega_{ij} \in (A_{s_{i}s_{j}})^{\times}$, mais il faut alors imposer les conditions $\omega_{ij} = \omega_{ik}.\omega_{kj}$ dans les anneaux  $A_{s_{i}s_{j}s_{k}}$. Le cas de deux \'el\'ements est nettement plus simple puisque ces conditions sont alors sans objet, et il met en sc\`ene d\'ej\`a  une partie de ce qui est en cause (pour le cas g\'en\'eral, voir {\sc Waterhouse} [9], 17.4).

\vspace{1cm}


{\Large{\bf 6. Application : le module inversible associ\'e \`a une forme quadratique binaire}}
\bigskip

Une forme quadratique binaire est une expression de la forme
$$
F(X, Y) = aX^2 + 2bXY + cY^2 . \leqno{(\star)}
$$
On lui associe l'\'el\'ement
$$
D = b^2-ac.
$$
L'\'etude syst\'ematique de ces formes remonte \`a Le Gendre et surtout \`a Gauss ; lorsque les coefficients sont des entiers, et que l'on consid\`ere comme \'equivalentes deux formes qui se d\'eduisent l'une de l'autre par un changement de variables lin\'eaire \`a coefficients entiers de d\'eterminant $+1$, alors Gauss a montr\'e qu'il n'y a qu'un nombre fini de formes non \'equivalentes. Ensuite Dedekind, \`a qui l'on doit la notion d'id\'eaux d'un corps de nombres, a montr\'e que la classification des formes quadratiques binaires \'etait essentiellement \'equivalente \`a celle des classes d'id\'eaux de l'anneau des entiers de $\Q[{\sqrt{D}}]$( voir {\sc Hecke}, [4] \S 53, Thm 154).\medskip 

Les remarques qui suivent \'etablissent une correspondance tr\`es naturelle entre forme quadratique binaire et module inversible sur cet anneau d'entiers alg\'ebriques ; cela clarifie le r\'esultat de Dedekind si on se souvient que le groupe des classes d'id\'eaux est isomorphe \`a celui des classes d'isomorphisme de modules inversibles.
\bigskip

\noindent ${\bf 6.1.}$\; L'id\'ee de  d\'epart, dans sa d\'esarmante simplicit\'e, consiste \`a lire la forme quadratique comme un d\'eterminant.

En effet, soit $R$ un anneau commutatif, et soit  $F(X, Y)$ une forme quadratique binaire \`a coefficients dans $R$. On a :
$$
F(X, Y) = X(aX+bY) -Y(-bX-cY) = \det \begin{pmatrix}
X &\, -bX-cY\\
Y&\, aX+bY
\end{pmatrix}.
$$
Introduisons l'endomorphisme $f$ du $R$-module $L = R^2$ d\'efini par sa matrice
$$
f =  \begin{pmatrix}
-b &\, -c\\
a&\, b
\end{pmatrix}.
$$
On a $\det(f) = ac - b^2 = - D$, et  ${\rm Tr}(f) = 0$ ; le th\'eor\`eme de Hamilton-Cayley montre que $f^2 = D$, et donc que $L$ peut \^etre muni d'une structure de module sur l'anneau 
$$
A = R[T]/(T^2 - D).
$$
On notera $\sqrt{D}$ la classe de $T$ dans $A$ ; pour $z \in L$, on a donc $\sqrt{D}.z = f(z)$.

Pour $ z = (x, y) \in R^2 = L$, l'application $A \rightarrow L,  \, \alpha \mapsto \alpha z$ a pour matrice, relativement \`a la $R$-base $\{1, \sqrt{D}\}$ de $A$, et \`a la base canonique de $L$,
$$
\begin{pmatrix}
x &\, -bx-cy\\
y&\, ax+by
\end{pmatrix}.
$$
Ainsi,
$$
F(x, y) \; = \; {\rm det}_{R}(A\,  \stackrel{\alpha \mapsto \alpha z}{\longrightarrow}\,  L).
$$
Par suite, un \'el\'ement $z = (x, y) \in L$ est une base du $A$-module $L$ si et seulement si $F(x, y)$ est inversible dans $R$. 
\bigskip

{\bf Proposition 6.2.}\; {\it Soient $a, b$ et $c$  des \'el\'ements d'un anneau $R$ tels que $aR+2bR+cR = R$. Alors le $A$-module $L$ introduit ci-dessus est inversible.}
\medskip

(Sous ces hypoth\`eses sur $a$, $b$ et $c$, Gauss nommait la forme \emph{proprement primitive}).
 
 Pour v\'erifier cette assertion, remarquons d'abord que si $a$, resp.  $c$, est inversible dans $R$ alors $(1, 0)$, resp. $(0, 1)$, est une base du $A$-module $L$ ; enfin, si $a+2b+c$ est inversible dans $R$ alors $(1, 1)$ est une base de $L$. Mais l'hypoth\`ese implique que $aR+(a+2b+c)R+cR = R$ ; il suffit donc d'invoquer {\bf 5.1.2.} pour pouvoir conclure.
 \bigskip
 
{\bf 6.3.}\; Montrons maintenant comment associer une forme quadratique binaire \`a un module inversible. Cette construction, qui utilise des rudiments d'alg\`ebre ext\'erieure, est un cas tr\`es particulier d'une construction g\'en\'erale d\'evelopp\'ee ailleurs\footnote{{\sc D. Ferrand}, {\it Un foncteur norme}, Bull. Soc. math. France, {\bf 126}, 1998, p. 1-49}.
\medskip

Soit $R \rightarrow A$ un morphisme d'anneaux faisant de $A$ un $R$-module libre de rang deux, et soit $L$ un $A$-module inversible ; alors, en consid\'erant les carr\'es ext\'erieurs des $R$-modules $A$ et $L$, on voit que $ \bigwedge^2 A$ est un $R$-module libre de rang $1$, et que ${\bigwedge^2} L$ est un $R$-module inversible, de sorte que 
$$
{\sf N}(L) = {\rm Hom}_{R}(\wedge^2 A, \wedge^2 L)
$$
 est un $R$-module inversible, d'ailleurs isomorphe, non canoniquement, \`a $\bigwedge^2L$. Par ailleurs, pour $z\in L$, le carr\'e ext\'erieur de l'application $A \rightarrow L,\; \alpha \mapsto \alpha.z$, est un \'el\'ement de ${\rm Hom}_{R}(\bigwedge^2 A, \bigwedge^2 L)$, que l'on note $\nu(z)$. Cela d\'efinit une application
$$
\nu : L \longrightarrow {\sf N}(L). \leqno{(\star\star)}
$$
C'est la forme quadratique binaire promise !
\bigskip

Avant de justifier cette affirmation p\'eremptoire, il faut souligner que l'application $\nu$ n'est  pas  additive, tout comme $F$, et que son image n'est en g\'en\'eral pas un $R$-module ; d'ailleurs,  c'est un probl\`eme redoutable de d\'eterminer l'ensemble ${\rm Im}(\nu)$, m\^eme dans les cas arithm\'etiques les plus simples, o\`u il  s'agit alors de caract\'eriser les entiers que l'on peut \'ecrire sous la forme $F(x, y)$ avec $x, y \in \Z$ (voir {\sc Landau} [6], Part Four, ch. IV).

Pla\c{c}ons-nous d'abord  dans la situation de {\bf 6.2.}, o\`u $L$ est le module inversible associ\'e \`a la forme $F$. La base canonique de $L = R^2$ conduit \`a un isomorphisme $\bigwedge^2L\;  \widetilde{\longrightarrow}\; R$, et la base $\{1, \sqrt{D}\}$ de $A$ donne un isomorphisme $\bigwedge^2A  \;\widetilde{\longrightarrow}\; R$ ; ces isomorphismes permettent d'identifier ${\sf N}(L)$ et  $R$. Il est alors clair que cela permet d'identifier les applications $\nu$ et $F$.

Pour ne pas imposer au lecteur l'\'elargissement de la notion de forme quadratique requis dans le cas g\'en\'eral, on va supposer que $L$ est libre sur $R$, donc de rang 2 ; le choix d'une base $(e_{1}, e_{2})$ de $L$ permet, comme ci-dessus, d'identifier les $R$-modules   ${\sf N}(L)$  et $R$ ; soit $\{1, t\}$ une base de $A$ comme $R$-module ; pour $z \in L$, l'\'el\'ement  $\nu(z) \in R$ est alors caract\'eris\'e par l'\'egalit\'e
$$
z \wedge t.z = \nu(z)e_{1}\wedge e_{2} .
$$
Introduisons la matrice de l'endomorphisme $z \mapsto t.z$, soit $\left( {\alpha \atop \beta} {\gamma \atop \delta} \right)$ ; on trouve
$$
\nu(Xe_{1}+Ye_{2}) =\beta X^2 + (\delta - \alpha)XY - \gamma Y^2
$$
Notons que si on peut choisir  le g\'en\'erateur $t$ de la $R$-alg\`ebre $A$ tel que ${\rm Tr}(t) = 0$, alors $\alpha + \delta = 0$, donc le coefficient $\delta - \alpha$ de $XY$  est bien un multiple de 2.
\bigskip

{\bf 6.4.}\;  Montrons, pour finir, comment d\'eterminer, dans le langage des modules inversibles, les automorphismes d'une forme quadratique binaire (cf {\sc Landau} [6], Thm 202, p. 181).

On suppose ici que $2$ est simplifiable dans l'anneau $R$. Soit, comme ci-dessus, $A$ une $R$-alg\`ebre libre de rang 2, et $L$ un $A$-module inversible ; on pr\'ecisera $\nu$ en $\nu_{L}$ pour l'application   $(\star\star)$, $L \rightarrow {\sf N}(L)$, introduite plus haut. 

Remarquons d'abord que ${\sf N}(A)$  est canoniquement isomorphe \`a $R$, et que l'application $\nu_{A}: A \rightarrow {\sf N}(A)=R$ n'est autre que l'application norme usuelle ; de sorte que l'\'equation
$$
\nu_{A}(\alpha) \; = \; 1
$$
est habituellement nomm\'ee  \emph{\'equation de Pell-Fermat } (relative au discriminant de $A/R$).

On v\'erifie imm\'ediatement que pour $u \in A$, on a
$$
\nu_{L}(uz) = \nu_{A}(u)\nu_{L}(z).
$$
Par suite, le produit dans $L$ par un $u \in A$ de norme $1$ induit un automorphisme de $L$ qui laisse $\nu_{L}$ invariante. La r\'eciproque est vraie :
\medskip

{\bf Proposition  6.5.}\; {\it Gardons les hypoth\`eses et les notations de {\bf 6.4.}. Soit $u$ un automorphisme $R$-lin\'eaire de $L$, de d\'eterminant $1$, et tel que $\nu_{L} \circ u = \nu_{L}$. Alors $u$ est une homoth\'etie de rapport un \'el\'ement de $A$ de norme $1$.}
\medskip

Il suffit de montrer qu'un automorphisme v\'erifiant ces deux conditions  est  $A$-lin\'eaire, puisque, $L$ \'etant un $A$-module inversible, les endomorphismes  $A$-lin\'eaires de $L$ sont les homoth\'eties.

Soit $t \in A$ tel que $\{1, t\}$ soit une base sur $R$. Il faut v\'erifier que, pour $z \in L$, on a $u(tz) = tu(z)$. Or, $1\wedge t$ est une base de $\bigwedge^2A$, et, par d\'efinition, $\nu_{L}(z)$ est l'application $\bigwedge^2A \rightarrow \bigwedge^2L, \; 1\wedge t \mapsto z\wedge t.z$ ; l'hypoth\`ese se traduit donc en : pour tout $z \in L$, on a 
$$
z\wedge t.z \; =\; u(z)\wedge t.u(z).
$$
Comme $u$ est suppos\'e de d\'eterminant $1$, soit $\wedge^2u = {\rm Id} = \wedge^2 u^{-1}$, on a  
$$
 z \wedge t.z =  z \wedge u^{-1}(t.u(z)),
$$
ce qui entra\^ine, pour tout $z \in L$,
$$
z \wedge (tz - u^{-1}(t.u(z))) \; =\; 0.
$$
Il reste \`a v\'erifier que le second facteur est nul.
\medskip

{\bf Lemme 6.6.}\; {\it Soit $R$ un anneau et  $f : L\rightarrow L$ un endomorphisme d'un $R$-module localement libre de rang $2$. Si, pour tout $z \in L$, on a $z \wedge f(z) = 0$, alors $f$ est une homoth\'etie de rapport $\lambda \in R$.}
\medskip

Donnons une d\'emonstration lorsque $L$ est libre, de base $\{e_{1}, e_{2}\}$ ; soit $\left( {\lambda\atop \mu} {\nu\atop \rho}\right)$ la matrice de $f$ sur cette base. Les relations $e_{i}\wedge f(e_{i}) = 0$ montrent que l'on a $\mu = \nu = 0$ ; en appliquant l'hypoth\`ese avec $z = e_{1}+e_{2}$, on obtient $e_{1}\wedge f(e_{2})+e_{2}\wedge f(e_{1}) = 0$, soit $\lambda = \rho$. $\Box$
\bigskip

Terminons la d\'emonstration de {\bf 6.5.} : d'apr\`es le lemme, il existe $\lambda \in R$ tel que $t.z - u^{-1}(t.u(z)) = \lambda z$ ; comme la trace des applications semblabes $z \mapsto u^{-1}(t.u(z))$ et $z\mapsto tz$  sont \'egales, on a $0 = {\rm Tr}(z\mapsto \lambda z) = 2\lambda$, d'o\`u $\lambda = 0$ puisque $2$ est suppos\'e simplifiable dans $R$.$\Box$
\vspace{3cm}

\noindent{\Large {\bf 7. Trivialit\'e des modules inversibles sur un anneau factoriel}}\\

\noindent {\bf Lemme 7.1.}\; {\it Soit $A$ un anneau int\`egre, et $L$ un $A$-module inversible. Alors :

\noindent {\rm 7.1.1.}\; Pour $a \in A$  et $x \in L$, la relation $ax = 0$ implique $a = 0$  ou $x = 0$.

\noindent {\rm 7.1.2.}\; Toute forme lin\'eaire non nulle $\alpha : L \rightarrow A$ est injective.}
\medskip

Pour v\'erifier ces propri\'et\'es, consid\'erons de nouveau des applications
$$
L\; \stackrel{u}{\longrightarrow}\; A^n\; \stackrel{v}{\longrightarrow}\; L \quad {\rm telles ~que} \quad v\circ u = {\rm Id}_{L}.
$$
La relation $ax = 0$ implique $au(x) = 0$, donc $a = 0$ ou $u(x) = 0$, puisque $A$ est int\`egre. Mais $x = v(u(x))$, d'o\`u (7.1.1).

Soit $\alpha : L \rightarrow A$ une forme lin\'eaire, et $x \in L$ tels que $\alpha(x) \neq 0$. Pour tout \'el\'ement non nul $y \in L$, on a, d'apr\`es  ce qui pr\'ec\`ede et 5.1.3, $0 \neq \alpha(x)y = \alpha(y)x$ ; donc $\alpha(y) \neq 0$.
\bigskip

De la factorialit\'e d'un anneau $A$, on n'utilisera que les propri\'et\'es suivantes :

- l'anneau $A$ est int\`egre ;

- pour tous $a, b \in A$, l'id\'eal  $aA \cap bA$ est principal.

\noindent En particulier, soit $\xi \in K$ un \'el\'ement non nul du corps des fractions de $A$. Alors, l'id\'eal 
$$\mathfrak{c} = \{ x \in A, x\xi \in A \}
$$
 est principal. En effet, si on \'ecrit $\xi = b/a$,  on a  $a.\mathfrak{c} = aA \cap bA$ ; un g\'en\'erateur de cet id\'eal principal est de la forme $ac$, avec $c \in \mathfrak{c}$ ; pour $x \in \mathfrak{c}$, il existe donc $y$ tel que $ax = acy$, d'o\`u $x = cy$ ; ainsi, $\mathfrak{c} = cA$.
\bigskip

{\bf Th\'eor\`eme 7.2.}\; {\it Si $A$ est un anneau factoriel, tout $A$-module inversible est isomorphe \`a $A$.}
\medskip

Soit $L$ un $A$-module inversible, et soit $\alpha : L \longrightarrow A$ une application lin\'eaire non nulle. Comme $\alpha$ est injective (7.1.2), cette application \'etablit un isomorphisme de $L$ sur l'id\'eal $\alpha(L)$. Il s'agit donc de montrer qu'un id\'eal $I \subset A$ qui est un $A$-module inversible, est principal.

Soit $\alpha : I \rightarrow A$ une forme lin\'eaire non nulle, et $x$ un \'el\'ement non nul de $I$. Pour tout $y \in I$, on a $\alpha(y)x = \alpha(yx) = y\alpha(x)$, soit
$$
\alpha(y) \; =\; \frac{\alpha(x)}{x}.y \; = \; \xi.y ,
$$
o\`u on a pos\'e $\xi = \alpha(x)/x$. 
Ainsi, avec un $x \in I$ fix\'e non nul, les formes lin\'eaires sur $I$ correspondent biunivoquement aux \'el\'ements $\xi$ du corps des fractions $K$ de $A$ tels que $\xi.I \subset A$.

Comme $I$ est projectif de type fini, il existe des applications
$$
I\; \stackrel{u}{\longrightarrow}\; A^n\; \stackrel{v}{\longrightarrow}\; I \quad {\rm telles ~que} \quad v\circ u = {\rm Id}_{I}.
$$
L'application $u$ est donn\'ee par $n$ formes lin\'eaires sur $I$ ; \`a chacune d'elle est associ\'ee, comme on vient de le voir, un \'el\'ement $\xi_{i} \in K$ tel que $\xi_{i}.I \subset A$.
Posons  $\mathfrak{c}_{i} = \{c \in A, c\xi_{i} \in A\}.$ Il a \'et\'e signal\'e plus haut que les $\mathfrak{c}_{i}$ sont des id\'eaux principaux, puisque $A$ est factoriel ; pour la m\^eme raison, une intersection finie d'id\'eaux principaux est un id\'eal principal ; il suffit donc de v\'erifier que 
$$
I\; =\; \mathfrak{c}_{1}\cap \cdots \cap \mathfrak{c}_{n} .
$$
Il est clair que $I$ est contenu dans cette intersection. R\'eciproquement, soit $a \in A$ tel que $a\xi_{i} \in A$ pour tout $i$. L'application $v : A^n \rightarrow I$ introduite plus haut, telle que $vu = {\rm Id}_{I}$, est d\'efinie par $n$ \'el\'ements $x_{1}, \ldots, x_{n}$ de $I$, et on a $\sum_{i} \xi_{i}x_{i} = 1$ ; par suite,
$
a \; =\; \sum_{i} (a\xi_{i})x_{i} \; \in \sum Ax_{i}\;  \subset \; I.
$
Cela ach\`eve la d\'emonstration.

\vspace{2cm}

\noindent {\Large{\bf 8. Modules inversibles et produit tensoriel}}
\bigskip

\noindent L'adjectif \emph{inversible} fait r\'ef\'erence au produit tensoriel, et est justifi\'e par la proposition suivante.
\bigskip

\noindent {\bf  Proposition  8.1.} {\it Soit $L$ un $A$-module pour lequel il existe un module $M$ et un isomorphisme}
$$
f : L\otimes_{A}M \quad \widetilde{\longrightarrow} \quad A.
$$
{\it Alors $L$ est inversible. R\'eciproquement, si $L$ est inversible, son dual ${\rm Hom}_{A}(L, A)$ est inversible, et l'application canonique}
$$
L\otimes_{A}{\rm Hom}_{A}(L, A) \; \longrightarrow \; A, \quad  x\otimes \alpha \longmapsto \alpha(x),
$$
{\it est un isomorphisme.}\medskip

Notons $g$ l'automorphisme de $L$ d\'efini par la commutativit\'e du diagramme suivant, o\`u on a d\'esign\'e par $\pi : M\otimes_{A}L\, \rightarrow\, L\otimes_{A}M$  l'isomorphisme de permutation.
$$
\begin{CD}
L\otimes_{A}M\otimes_{A}L @>{1\otimes \pi}>> L\otimes_{A}L\otimes_{A}M\\
@V{f\otimes 1}VV  @VV{1\otimes f}V\\
L @>>g> L\\ 
\end{CD}
$$
Pour $x, y \in L$  et $z\in M$, en suivant les destins de $y\otimes z \otimes x \in L\otimes_{A}M\otimes_{A}L$,  on voit que l'on a
$$
f(y\otimes z)g(x)\; = \; yf(x\otimes z).
$$
Par ailleurs,la surjectivit\'e de $f$ entra\^ine  l'existence d'\'el\'ements $y_{1}, \ldots , y_{n}  \in L$, et d'\'el\'ements $z_{1}, \ldots , z_{n} \in M$ tels que
$$
\sum f(y_{i}\otimes z_{i}) \, = \, 1.
$$

De ces deux \'egalit\'es, on tire que, pour tout $x \in L$, on a 
$$
g(x) = (\sum f(y_{i}\otimes z_{i}))g(x) = \sum y_{i}f(x\otimes z_{i}). \leqno{(\star)}
$$
Or, les $n$ formes lin\'eaires $x \mapsto f(x\otimes z_{i})$ d\'efinissent une application $v : L \rightarrow A^n$, et les \'el\'ements $y_{i}$ d\'efinissent une application $u : A^n \rightarrow L$. L'\'egalit\'e $(\star)$ s'\'ecrit alors
$$
g \, = \, u\circ v .
$$
Comme $g$ est un isomorphisme, cela montre que $L$ est facteur direct de $A^n$.

Enfin, comme la dimension d'un produit tensoriel d'espaces vectoriels est \'egale au produit des dimensions des facteurs, il est clair que $L$ est de rang 1.

R\'eciproquement, soit $L$ un module inversible ; notons $L\check{{}}$\, son dual. Introduisons, comme en 4.1.1 des applications lin\'eaires
$$
L \stackrel{u}{\longrightarrow}  A^n  \stackrel{v}{\longrightarrow}  L, \qquad {\rm telles ~ que}\; v\circ u = {\rm id}_{L}.
$$
Passant aux duals, on trouve des applications lin\'eaires
$$
L\check{{}}\,  \stackrel{u\check{{}}}{\longleftarrow}  \, (A^n)\check{{}} \, \stackrel{v\check{{}}}{\longleftarrow} \, L\check{{}}
$$
dont le compos\'e est l'identit\'e. Cela montre d\'ej\`a que $L\check{{}}$ est projectif de type fini ; le rang de ce module est \'egal d'apr\`es 4.1.3 au rang, en chaque id\'eal premier, de l'application compos\'ee
$$
(A^n)\check{{}} \, \stackrel{v\check{{}}}{\longleftarrow} \, L\check{{}} \stackrel{u\check{{}}}{\longleftarrow}  \, (A^n)
$$
Mais c'est la duale de l'application $uv$, laquelle est de rang 1. Donc $L\check{{}}$\, est inversible.

Il en resulte que $L\otimes_{A}L\check{{}}$ est inversible, si bien qu'il suffit, d'apr\`es 4.1.7,  de montrer la surjectivit\'e de l'application $L\otimes_{A}L\check{{}} \, \rightarrow \, A$ pour pouvoir conclure que c'est un isomorphisme. Or, la relation $v\circ u = {\rm id}_{L}$ se traduit ainsi : il existe $n$ formes lin\'eaires $u_{i} : L \rightarrow A$ et $n$ \'el\'ements $x_{i} \in L$ tels que, pour tout $x \in L$, on a
$$
\sum_{i}u_{i}(x)x_{i} \, =\, x .
$$
D'apr\`es 5.1.3, on a $u_{i}(x)x_{i} = u_{i}(x_{i})x$, d'o\`u
$$
(\sum_{i} u_{i}(x_{i})).x = x
$$
On en tire que $\sum_{i} u_{i}(x_{i}) = 1$, donc que l'image de $\sum_{i}x_{i}\otimes u_{i} \in L\otimes_{A}L\check{{}}$ par l'application en cause est \'egale \`a 1.
$\Box$ 
\bigskip

\noindent {\bf 8.2.}\; Les classes d'isomorphie de $A$-modules inversibles forment donc un groupe pour le produit tensoriel ; l'\'el\'ement neutre en est la classe de $A$, et l'inverse de la classe de $L$ est la classe de son dual.

Ce groupe est nomm\'e le groupe de Picard de $A$, et not\'e 
$${\rm Pic}(A).$$
 Il joue un r\^ole consid\'erable, tant en g\'eom\'etrie alg\'ebrique qu'en th\'eorie alg\'ebrique des nombres o\`u il est plut\^ot vu comme groupe des classes d'id\'eaux.

\vspace{1cm}

{\Large {\bf 9. Descente}}
\bigskip

La compr\'ehension de ce paragraphe requiert tr\`es peu de connaissances sur les morphismes fid\`element plats ; pour en lire plus, on consultera le livre de {\sc Knus-Ojanguren} [5], ou le ch. 17 de celui de {\sc Waterhouse} [9].

On utilisera essentiellement la caract\'erisation suivante : pour qu'un morphisme d'anneaux  $A \rightarrow B$ soit fid\`element plat, il faut et il suffit que la propri\'et\'e suivante soit v\'erifi\'ee :
\medskip

\noindent {\bf 9.1.}\;  {\it  Pour qu'une application $A$-lin\'eaire $ E\rightarrow F$ soit injective (resp. surjective, resp. un isomorphisme), il faut et il suffit que l'application $B$-lin\'eaire $B\otimes_{A}E \, \rightarrow \, B\otimes_{A}F$ soit injective (resp. \ldots ).}
\medskip

On en d\'eduit l'important r\'esultat suivant :
\medskip

\noindent {\bf 9.2.}\; {\it Soit $ u : A \rightarrow B$ un morphisme fid\`element plat. Alors, le morphisme $u$ est injectif et $u(A)$ est \'egal \`a l'ensemble des \'el\'ements $b \in B$ tels que, dans $B\otimes_{A}B$, on ait  $ b\otimes 1 = 1 \otimes b$. Ce qu'on r\'esume en disant que la suite de morphismes d'anneaux ci-dessous est exacte, o\`u on a not\'e $u_{0}$ l'application $b \mapsto 1\otimes b$, et $u_{1}$ l'application $b \mapsto b\otimes 1$ :}
$$
\xymatrix{
A \ar[r]^u & B\ar@< 2pt>[r]^-{u_0} \ar@<-3pt>[r]_-{u_1} & B\otimes_{A}B
}
$$
\medskip

Dans le contexte des $A$-modules on \'ecrirait plut\^ot que la suite
$$
0 \; \longrightarrow \; A \;\stackrel{u}{\longrightarrow}\; B \; \stackrel{u_{1}-u_{0}}{\longrightarrow}\; B\otimes_{A}B
$$
est exacte ; mais cela risque de faire oublier que $u,\, u_{0}$  et $u_{1}$ sont des morphismes d'alg\`ebres, ce que $u_{1} - u_{0}$ n'est pas.

Indiquons la d\'emonstration de ce r\'esultat, d\^u \`a Grothendieck, et qui devrait figurer dans tous les manuels d'alg\`ebre commutative depuis quarante ans (voir aussi {\sc Knus-Ojanguren} [5], p.30).
\medskip

D'apr\`es (9.1) il suffit de v\'erifier que cette suite de morphismes devient exacte apr\`es tensorisation \`a droite par $B$ , soit l'exactitude de
$$
\xymatrix{
B \ar[r]^{u \otimes 1} & B\otimes_{A}B \ar@< 2pt>[r]^-{u_0\otimes 1} \ar@<-3pt>[r]_-{u_1\otimes 1} & B\otimes_{A}B\otimes_{A}B
}
$$
Dans ce diagramme, et dans les suivants, le symbole $1$ d\'esignera souvent l'application identique d'un module que le contexte d\'esignera clairement ; pr\'ecisons cependant que le facteur $A\otimes_{A}-$ est omis, si bien que $u\otimes 1$ d\'esigne le morphisme compos\'e
$$
B \stackrel{b\mapsto 1\otimes b}{\longrightarrow} A\otimes_{A}B \; \stackrel{u\otimes 1}{\longrightarrow}\, B\otimes_{A}B.
$$
On a donc $(u\otimes 1)(b) = 1\otimes b = u_{0}(b)$.
 D\'esignons par $s : B\otimes_{A}B \rightarrow B$ le morphisme donn\'e par le produit : $x\otimes y \mapsto xy$, et d\'esignons par $t = 1\otimes s : B\otimes_{A}B\otimes_{A}B \rightarrow B\otimes_{A}B$, celui donn\'e par $x\otimes y\otimes z \mapsto x\otimes yz$. On constate que $s \circ (u\otimes 1) = {\rm Id}_{B}$, donc que les morphismes  $u\otimes 1$, et par suite aussi  $u$ sont injectifs ; on constate aussi que 
$$
t \circ (u_{0}\otimes 1) = u_{0}\circ s, \qquad {\rm et}\quad t \circ (u_{1}\otimes 1) = {\rm Id}_{B\otimes_{A}B}
$$
Si, donc, un \'el\'ement $\xi \in B\otimes_{A}B$ v\'erifie $(u_{0}\otimes 1)(\xi) = (u_{1}\otimes 1)(\xi)$, en composant avec $t$, on trouve
$$
u_{0}(s(\xi)) = 1\otimes s(\xi) = \xi
$$
C'est ce qu'on voulait v\'erifier.$\Box$
\vspace{1cm}

{\bf Th\'eor\`eme 9.3.}\;  {\it Soit $A \rightarrow B$ un morphisme fid\`element plat d'anneaux, et $L$ un $A$-module. S'il existe un isomorphisme de $B$-modules $B\otimes_{A}L\;  \widetilde{\longrightarrow}\;  B$, alors $L$ est un $A$-module inversible.}
\bigskip

La suite exacte de 9.2 conduit \`a la suite exacte 
$$
{\rm Hom}_{A}(L, A)\; \rightarrow\; {\rm Hom}_{A}(L, B)\; \rightrightarrows \; {\rm Hom}_{A}(L, B\otimes_{A}B) .
$$

La  propri\'et\'e (9.1) implique visiblement que $L$ est un $A$-module plat, si bien que, par tensorisation \`a droite par $L$,  on obtient la suite exacte
$$
{\rm Hom}_{A}(L, A)\otimes_{A} L\; \rightarrow\; {\rm Hom}_{A}(L, B)\otimes_{A} L\; \rightrightarrows \; {\rm Hom}_{A}(L, B\otimes_{A}B)\otimes_{A} L .
$$
On va comparer cette suite \`a celle de 9.2, en utilisant le r\'esultat suivant :\\

\noindent {\bf 9.4.}\; {\it Soit $A\, \rightarrow\, C$ une $A$-alg\`ebre pour laquelle on dispose d'un isomorphisme de $C$-modules 
$$\omega : C\otimes_{A}L\, \tilde{\rightarrow}\, C.$$
 Alors, l'application $w_{C} : {\rm Hom}_{A}(L, C)\otimes_{A}L\, \rightarrow \, C$, d\'efinie par $ \alpha\otimes x \mapsto \alpha(x)$ est un isomorphisme.}
\medskip

Consid\'erons, en effet,  la suite d'isomorphismes
$$
C \, \simeq \, {\rm Hom}_{C}(C, C)\, \stackrel{\circ \omega }{\longrightarrow}\, {\rm Hom}_{C}(C\otimes_{A} L, C)\, \stackrel{\varphi}{\longrightarrow}\, {\rm Hom}_{A}(L, C).
$$
Le symbole $\circ \omega$ d\'esigne la composition \`a droite par l'isomorphisme $\omega$, et $\varphi$ d\'esigne l'isomorphisme canonique d\'eduit de la composition avec l'application $L \rightarrow C\otimes_{A} L$. L'isomorphisme compos\'e $\psi : C \rightarrow {\rm Hom}_{A}(L, C)$ est celui qui \`a $c \in C$ associe l'application $A$-lin\'eaire  d\'efinie par $x \mapsto \omega(c\otimes x)$.

Le carr\'e suivant est commutatif :
$$
\begin{CD}
C\otimes_{A}L @>{\omega}>> C\\
@V{\psi\otimes 1}VV    @|\\
{\rm Hom}_{A}(L, C)\otimes_{A}L @>>{w_{C}}>C
\end{CD}
$$
Par suite, l'application du bas est un isomorphisme.$\Box$
\medskip

Consid\'erons alors le diagramme commutatif suivant :
$$
\xymatrix{
{\rm Hom}_{A}(L, A)\otimes_{A}L \ar[r]^{u'}\ar[d]_{w_{A}} & {\rm Hom}_{A}(L, B)\otimes_{A}L \ar[d]_{w_{B}} \ar@< 2pt>[r]^-{u'_0} \ar@<-3pt>[r]_-{u'_1}
& {\rm Hom}_{A}(L, B\otimes_{A}B)\otimes_{A}L
 \ar[d]_{w_{B\otimes_{A}B}} \\
A \ar[r]_u & B\ar@< 2pt>[r]^-{u_0} \ar@<-3pt>[r]_-{u_1} & B\otimes_{A}B
}
$$

Les lignes sont exactes et  $w_{B}$ et $w_{B\otimes_{A}B}$ sont des isomorphismes d'apr\`es le r\'esultat pr\'ec\'edent appliqu\'e avec $C = B$ et $C = B\otimes_{A}B$. Par suite, $w_{A}$ est un isomorphisme, et $L$ est un module inversible (8.1).

\vspace{1cm}


{\Large {\bf  10. Constructions galoisiennes}}
\bigskip

\noindent {\bf 10.1.}\;  Le th\'eor\`eme pr\'ec\'edent peut \^etre prolong\'e en un proc\'ed\'e de \emph{construction} de modules inversibles, \`a partir d'un morphisme fid\`element plat $A \rightarrow B$, et de ce qu'on appelle une \emph{donn\'ee de descente} ; on n'abordera pas ici cette construction sous sa plus grande g\'en\'eralit\'e (voir {\sc Knus-Ojanguren} [5], p.36, ou {\sc Waterhouse} [9], p.132) ; on se limitera \`a des morphisme tr\`es particuliers, les rev\^etements galoisiens, parce que la \og donn\'ee de descente\fg\,  se r\'eduit alors \`a un cocycle (voir aussi {\sc Knus-Ojanguren} [5], p.44, ou {\sc Waterhouse} [9], p.136).
\medskip

Soit $B$ un anneau (commutatif) muni d'un groupe fini $G$ d'automorphismes, et soit $ A = B^G$ le sous-anneau des \'el\'ements invariants. On cherche \`a construire des  $A$-modules $L$ munis d'un isomorphisme  de $B$-modules 
$$
\omega : B\otimes_{A}L \; \widetilde{\longrightarrow}\; B.
$$
Supposons donn\'e un tel isomorphisme. Un automorphisme $g \in G$ est $A$-lin\'eaire par d\'efinition de $A$ ; il induit, par suite, un isomorphisme $A$-lin\'eaire $g\otimes 1 : B\otimes_{A}L \longrightarrow B\otimes_{A}L$. On d\'efinit une application $A$-lin\'eaire $\varphi$ (qui d\'epend de $g$) par la commutativit\'e du carr\'e suivant :
$$
\begin{CD}
B\otimes_{A}L @> g\otimes 1 >> B\otimes_{A}L\\
@V \omega VV @VV\omega V\\
B @>>\varphi > B\\
\end{CD}
$$
Pour $b, c \in B$, on a $\varphi(bc) = g(b)\varphi(c)$, puisque cela est vrai pour l'application $g \otimes 1$, et que $\omega$ est un isomorphisme $B$-lin\'eaire.  Posons $\varphi(1) = \theta(g)$, de sorte que $\varphi(b) = \theta(g).g(b).$ On v\'erifie imm\'ediatement la propri\'et\'e suivante :
$$
\theta(1) = 1,\; {\rm et\, pour\, tout } \; g, h \in G,  \; {\rm on\, a} \quad  \theta(gh) = \theta(g).g(\theta(h)). \leqno{({\sf C})}
$$
Elle implique que les $\theta(g)$ sont inversibles dans $B$. 
\medskip

\noindent Une application $\theta : G \rightarrow B^{\times}$ satisfaisant les relations ({\sf C}) est nomm\'ee un \emph{cocycle} de $G$ \`a valeurs dans $B$. Si $\theta$ est \`a valeurs dans $A$, c'est-\`a-dire si les $\theta(g)$ sont invariants sous $G$, la condition {\sf C} signifie simplement que  $\theta$ est un homomorphisme de groupes $G \rightarrow A^{\times}$. Notons aussi que, pour un \'el\'ement inversible $u \in B^{\times}$, l'application d\'efinie par $\theta(g) = u/g(u)$ est un cocycle ; ces cocycles sont nomm\'es des \emph{cobords}.
\medskip

Revenons au module $L$. L'application compos\'ee $L \stackrel{x \mapsto 1\otimes x}{\longrightarrow} B\otimes_{A}L \stackrel{\omega}{\longrightarrow} B$ donne  une application  
$$
L \, \longrightarrow \,\{b\in B, \forall g\in G \quad   \theta(g).g(b) = b\}.
$$
 L'id\'ee directrice du paragraphe est de d\'efinir/construire le module \`a partir d'un cocycle $\theta$, en \emph{posant}
 $$
L_{\theta}  = \{b\in B, \forall g\in G \quad   \theta(g).g(b) = b\},\leqno{10.1.1}
$$
et de d\'egager des conditions sur le morphisme $A\rightarrow B$ qui assureront que le $A$-module $L_{\theta}$ ainsi construit est inversible.

\medskip

{\bf 10.2.}\; Pour un ensemble fini $J$, on note $\prod_{J}B$ l'anneau produit de copies de $B$ index\'ees par les \'el\'ements de $J$ ; un \'el\'ement de cet anneau est donc une famille $(b_{j})_{j\in J}$ d'\'el\'ements de $B$, qu'on d\'ecrira le plus souvent  comme une application $J\rightarrow B$ (il est clair que la notation fr\'equente $B^J$ entre en conflit avec les invariants).

\noindent \`A l'action d'un groupe fini $G$ sur $B$, et en posant $A = B^G$, on associe le morphisme d'anneaux
$$
\rho_{B, G} : B\otimes_{A}B\; \longrightarrow \; \prod_{G} B,\qquad  x\otimes y \mapsto (g\mapsto xg(y)).
$$

\bigskip

{\bf D\'efinition 10.3.} {\it Soit $G$ un groupe fini. On dira qu'un morphisme $u : A\rightarrow B$ est \emph{galoisien de groupe} $G$ si ce groupe op\`ere sur $B$ de telle sorte que les conditions suivantes sont satisfaites}

1) {\it on a $B^G = A$ ;}

2) {\it le morphisme $\rho_{B, G}$ est un isomorphisme} ;

3)  {\it le morphisme $u$ fait de $B$  un $A$-module projectif de type fini, n\'ecessairement de rang $d = {\rm Card}(G)$.}
\bigskip

Pour tout groupe fini $G$, et tout anneau $A$, la $A$-alg\`ebre $B = \prod_{G}A$ est donc un rev\^etement galoisien pour l'op\'eration de $G$ sur $B$ donn\'ee par : $(gb)(g') = b(g'g)$.
\medskip

{\bf Exercices 10.4.}\; a) Montrer qu'une extension finie galoisienne de corps $K \subset K'$, de groupe de Galois $G$,  est un morphisme galoisien de groupe $G$ (Pour montrer que le morphisme $\rho_{K', G}$ est un isomorphisme, on peut utiliser le th\'eor\`eme de l'\'el\'ement primitif : $K' = K(t)$, et le fait que $K'/K$ est une extension de d\'ecomposition du polyn\^ome minimal de $t$).

b) Soit $u : A \rightarrow B$ un morphisme v\'erifiant les propri\'et\'es 2) et 3) de la d\'efinition. En reprenant les notations et r\'esultats de 9.2, calculer $\rho \circ u_{0}$  et  $\rho \circ u_{1}$. En d\'eduire que la propri\'et\'e 1) est une cons\'equence  de 2) et 3).

c) Soit $u : A \rightarrow B$ un morphisme galoisien de groupe $G$. Montrer que pour tout $t \in A$, le morphisme $A/tA \rightarrow  B/u(t)B$ est encore galoisien de groupe $G$. Plus g\'en\'eralement, pour toute $A$-alg\`ebre $A\rightarrow A'$ le morphisme $A' \rightarrow A'\otimes_{A}B$ est galoisien de groupe $G$ (On pourra utiliser la question pr\'ec\'edente).
\bigskip

{\bf Proposition 10.5.}\; {\it Soient $u : A \rightarrow B$ un morphisme galoisien de groupe $G$, et \,$\theta : G \rightarrow B^{\times}$ un cocycle. Posons :}
$$
L_{\theta}  = \{b\in B, \forall g\in G \quad   \theta(g).g(b) = b\} .
$$
{\it Alors le $A$-module associ\'e  $L_{\theta}$ est inversible. Il est libre si et seulement si il contient un \'el\'ement qui est inversible dans $B$, autrement dit si $\theta$ est un } cobord.
\medskip

Montrons d'abord que le morphisme canonique
$$
\omega : B\otimes_{A}L_{\theta}\; \rightarrow \; B
$$
est un isomorphisme ; d'apr\`es le th\'eor\`eme de descente (9.3), cela impliquera que ce module est inversible.
\bigskip

Posons $L = L_{\theta}$. Le module  $L$ s'ins\`ere dans la suite exacte suivante de $A$-modules, o\`u on a not\'e $\psi$ l'application d\'efinie par $\psi(b) = (g\mapsto \theta(g).g(b) - b)$
$$
0\; \longrightarrow \; L \; \stackrel{\iota}{\longrightarrow}\; B \; \stackrel{\psi}{\longrightarrow}\; \prod_{G} B
$$
Consid\'erons le diagramme suivant d'applications $B$-lin\'eaires, o\`u $\rho = \rho_{B, G}$ est un isomorphisme par hypoth\`ese :
$$
\begin{CD}
0 @>>> B\otimes_{A}L @>{1\otimes \iota}>> B\otimes_{A}B @>{1\otimes \psi}>>  B\otimes_{A} \prod_{G}B\\
&& @V\omega VV @VV\rho V &\\
0 @>>> B @>>{1 \mapsto (g\mapsto \theta(g)^{-1})}> \prod_{G}B \\
\end{CD} \leqno 10.5.1
$$
Le carr\'e est construit pour \^etre commutatif, et il l'est ! Comme $B$ est un $A$-module libre de rang fini, la suite du haut  est exacte. Cela montre d\'ej\`a que $\omega$ est injectif puisque $\rho$ et $1\otimes   \iota$ le sont. Mais $\rho$ est aussi surjectif ; il existe donc un \'el\'ement $z  \in B\otimes_{A}B$ tel que $\rho(z)= (g\mapsto \theta(g)^{-1})$. On va v\'erifier que cet \'el\'ement $z$ est dans le noyau de $1\otimes \psi$, donc qu'il est dans $B\otimes_{A}L$ ; la commutativit\'e du carr\'e impliquera que l'on a   $\omega (z) = 1$, et cela d\'emontrera que $\omega$ est un isomorphisme. 
\medskip

En explicitant $z = \sum x_{i}\otimes y_{i}$, la propri\'et\'e $\rho(z) = (g\mapsto \theta(g)^{-1})$ se traduit en ceci : pour tout $g\in G$,
$$
\sum_{i} x_{i}\theta(g)g(y_{i}) \; =\; 1. \leqno{(\star)}
$$

\noindent Par d\'efinition, $(1\otimes \psi) (\sum_{i} x_{i}\otimes y_{i})$ est l'application de $G$ dans $B\otimes_{A}B$ d\'efinie par
$$
g \; \longmapsto\; \sum_{i} x_{i}\otimes (\theta(g)g(y_{i})-y_{i}).
$$
Notons $z(g)$ le membre de droite ; il s'agit de montrer que pour tout $g \in G$,  cet  \'el\'ement   $z(g) \in B\otimes_{A}B$ est nul. Fixons un \'el\'ement $g \in G$. Comme $\rho$ est injectif, il suffit de voir que $\rho(z(g))$ est nul. Or, on a
$$
\rho(z(g)) \;=\; \rho[\sum_{i} x_{i}\otimes (\theta(g)g(y_{i})-y_{i})] \; = (h \, \longmapsto \, \sum_{i} x_{i} h[\theta(g)g(y_{i})-y_{i}])
$$
Compte-tenu de la propri\'et\'e de cocyle {\sf (C)}, on constate que 
$$
 \sum_{i} x_{i} h[\theta(g)g(y_{i})-y_{i}] = \theta(h)^{-1}. \Big[ \sum_{i} x_{i} \theta(hg)hg(y_{i})\, - \,  \sum_{i} x_{i} \theta(h).h(y_{i})\Big]
 $$
 La nullit\'e de cet \'el\'ement provient des relations $(\star)$. 
 
 On a donc d\'emontr\'e que l'application $\omega : B\otimes_{A}L_{\theta}\; \rightarrow \; B$ est un isomorphisme.
 \medskip
 
 Supposons que $L$ contienne un \'el\'ement $u$ inversible dans $B$. Pour tout $b \in L$, on a, pour tout $g \in G$, \`a la fois $\theta(g).g(b) = b$  et $\theta(g).g(u) = u$, d'o\`u $g(b/u) = b/u$, c'est-\`a-dire $b \in Au$, puisque $B^G
 = A$ ; ainsi, $L = Au$.
 
 R\'eciproquement, supposons que $L$ soit libre, engendr\'e par $u$. La surjectivit\'e de $\omega$ implique qu'il existe des \'el\'ements $b_{i} \in B$ et des \'el\'ements $x_{i} \in L$, tels que $\sum_{i} b_{i}x_{i} = 1$ ; par hypoth\`ese, chaque $x_{i}$ est de la forme $a_{i}u$, avec $a_{i} \in A$ ; on a donc $(\sum_{i}b_{i}a_{i}).u = 1$, et $u$ est inversible dans $B$.$\Box$
 \bigskip
 
 {\bf Exemple 10.6.}\; Reprenons l'exemple d'une extension galoisienne finie de corps, $K \subset K'$, comme en 10.4. a), de sorte qu'ici $G = {\rm Gal}(K'/K)$. Un cocycle $\theta : G \rightarrow K'^{\times}$ d\'etermine un $K$-module inversible $L_{\theta}$, c'est-\`a-dire un espace vectoriel de dimension $1$, lequel est libre ! La proposition 10.5 montre donc que pour tout cocycle $\theta$, il existe un \'el\'ement  $u \in K'^{\times}$ tel que, pour tout $g \in G$,  
 $$
 \theta(g) \; =\; u/g(u) .
 $$
 Ce r\'esultat :  \emph{pour une extension de corps, tout cocycle est un cobord}, se d\'emontre facilement \`a l'aide d'une r\'esolvante de Lagrange ; il est quelquefois nomm\'e le \og th\'eor\`eme 90 de Hilbert.\fg
 \medskip
 
 {\bf Exemple 10.7.}\; Revenons \`a la situation du d\'ebut (2.2 et 3.5). L'anneau
 $$A = \R[X, Y]/(X^2+Y^2-1)
 $$
 est consid\'er\'e comme sous-anneau de 
 $$B = \C[X, Y]/(X^2+Y^2-1)$$ 
 Il est clair que $B = \C\otimes_{\R}A$, si bien que le morphisme $A \rightarrow B$ est galoisien de groupe $G = \{{\rm Id}, \sigma\}$, o\`u $\sigma$ est la conjugaison complexe.
 
 On d\'efinit un cocycle 
 $$
 \theta : G \; \longrightarrow \; B^{\times}
 $$
 en posant \hspace{4cm} $\theta({\rm Id}) = 1$, \quad  et \quad $\theta(\sigma) = x+iy$.
 
 C'est effectivement un cocycle puisque  $\theta(\sigma).\sigma(\theta(\sigma)) = (x+iy).(x-iy) = x^2+y^2 = 1$, et que $\theta({\rm Id}) = \theta(\sigma^2)$. La relation (2.2.1) dit exactement que le module $M$ associ\'e au ruban de M\"obius est form\'e des \'el\'ements $a+ib \in B$ tels que
 $$
 (x+iy).(a-ib) = a+ib, \qquad {\rm soit}\qquad  \theta(\sigma).\sigma(a+ib) = a+ib.
 $$
 C'est donc le module associ\'e au cocycle $\theta$.
 \bigskip
 
 {\bf (Pseudo) Exemple 10.8.}\; ({\it Signature})\, Posons $V(X_{1}, \ldots, X_{n}) = \prod_{i < j} (X_{j} - X_{i})$. Par d\'efinition de la signature $\varepsilon (\sigma)$ d'une permutation $\sigma \in \mathfrak{S}_{n}$, on a 
 $$
 V(X_{\sigma(1)}, \ldots, X_{\sigma(n)})\; = \; \varepsilon(\sigma).V(X_{1}, \ldots, X_{n}).
 $$
 Soit $P(X_{1}, \ldots, X_{n})$ un polyn\^ome \`a coefficients dans un corps $K$ de caract\'eristique $\neq 2$, tel que pour tout $\sigma \in \mathfrak{S}_{n}$, on ait
 $$
 P(X_{\sigma(1)}, \ldots, X_{\sigma(n)})\, =\, \varepsilon(\sigma).P(X_{1}, \ldots, X_{n}).
 $$
 Rappelons comment on v\'erifie que cette propri\'et\'e \'equivaut \`a : $P = V.Q$, o\`u $Q$ est un polyn\^ome sym\'etrique. Notons ${}^{\sigma}P$ le polyn\^ome $P(X_{\sigma(1)}, \ldots, X_{\sigma(n)})$, et all\`egeons les notations en d\'esignant par $F(X_{i}, X_{j})$ le polyn\^ome $P$ vu comme polyn\^ome en $X_{j}$ et $X_{j}$ \`a coefficients des polyn\^omes en les autres ind\'etermin\'ees. Pour la transposition $\tau = (i\, j)$, on a ${}^{\tau}P = -P$, donc $F(X_{j}, X_{i}) = - F(X_{i}, X_{j})$ ; faisant $X_{j} = X_{i}$, on obtient  $2F(X_{i}, X_{i}) = 0$, soit $F(X_{i}, X_{i}) = 0$ puisque $2$ est suppos\'e inversible dans $K$ ; cela montre que $P$ est multiple de $X_{j}-X_{i}$. Comme, pour deux couples distincts $(i, j)$  et $(i', j')$ ces polyn\^omes n'ont pas de diviseurs communs, on voit que $V$ divise $P$, soit $P = V.Q$ ; il est clair que $Q$ est sym\'etrique.
 \medskip
 
 Interpr\'etons ce r\'esultat fort classique dans le contexte du paragraphe. Posons $B = K[X_{1}, \ldots, X_{n}]$, et $A = B^{\mathfrak{S}_{n}} = K[S_{1}, \ldots, S_{n}]$, o\`u les $S_{i}$ sont les polyn\^ omes sym\'etriques \'el\'ementaires en les $X_{j}$. Le morphisme $A\rightarrow B$ n'est pas galoisien car la propriet\'e 2) de la d\'efinition n'est pas v\'erifi\'ee, mais la d\'emarche garde un sens (on peut v\'erifier que le morphisme localis\'e $A_{V} \longrightarrow B_{v}$ est galoisien). La signature $\varepsilon : \mathfrak{S}_{n} \rightarrow \{\pm 1\}$ est un cocyle, auquel est associ\'e le $A$-module 
 $$
 L_{\varepsilon} = \{P \in B, \forall \sigma \in \mathfrak{S}_{n},  \varepsilon(\sigma).{}^{\sigma}P = P \}
 $$
On vient de rappeler que $L_{\varepsilon} = V.A $. Cela confirme que c'est un $A$-module inversible, et qu'il est libre comme il se doit, puisque $A$ est factoriel.

 \vspace{5cm}
 
 {\Large {\bf  11. Application : la th\'eorie de Kummer}}
\bigskip

Elle d\'ecrit les rev\^etements galoisiens dont le groupe est ab\'elien.
\bigskip

{\bf 11.1.}\; Avant d'\'enoncer le r\'esultat en vue, donnons un exemple simple \`a mi-chemin entre la th\'eorie de Kummer classique et ce qui va suivre.

Soit $A$ un anneau et $L$ un $A$-module inversible ; on suppose donn\'es un entier $d$ et un isomorphisme
$$
\pi : L^{\otimes d} \; \widetilde{\longrightarrow}\; A.
$$
Pour un entier $k$ tel que $0 \leq k\leq d-1$, on d\'esigne par $\pi_{k}$ l'isomorphisme compos\'e
$$
L^{\otimes d+k} \simeq L^{\otimes d}\otimes L^{\otimes k}\; \stackrel{\pi \otimes 1}{\longrightarrow}\; L^{\otimes k}.
$$
Le choix du facteur $L^{\otimes d}$ sur lequel on applique $\pi$ est indiff\'erent, cf 5.1.3.

Posons 
$$
B = A \oplus L \oplus L^{\otimes 2} \oplus \cdots \oplus L^{\otimes d-1}.
$$
 Ce $A$-module est muni d'un produit d\'efini de la fa\c{c}on suivante : pour $x \in L^{\otimes i}$  et $y \in L^{\otimes j}$, on pose $x.y = x\otimes y \in L^{\otimes i+j}$, si $i+j < d$, et si $i+j \geq d$, on pose  $x.y = \pi_{i+j-d}(x\otimes y)$. Le lecteur (ne) v\'erifiera (pas) que cela fait de $B$ une $A$-alg\`ebre commutative. 

Lorsque $L$ est libre, de g\'en\'erateur $x$, alors l'\'el\'ement $t = \pi(x^{\otimes d})$ est inversible dans  $A$ puisque $\pi$ est un isomorphisme, et $B$ est visiblement isomorphe \`a $A[X]/(X^d - t)$, autrement dit, $B$ est alors une extension de Kummer au sens habituel. On verra plus bas que sous des condition assez g\'en\'erales $B$ est un rev\^etement galoisien de $A$.

\bigskip

\noindent {\bf 11.2.} \; D\'egageons d'abord les hypoth\`eses n\'ecessaires \`a la validit\'e de ce qui suit.
Je remercie L. Moret-Bailly de m'avoir fait remarquer qu'une pr\'ec\'edente version \'etait beaucoup trop optimiste.
\medskip

\noindent Soit $N$ un entier $ > 0$. D\'esignons par $({\sf P}_{N})$ la propri\'et\'e suivante qui porte sur un anneau commutatif~$A$~:
\medskip

\noindent $({\sf P}_{N})$ :{\it  Le groupe mutiplicatif $A^{\times}$ contient un sous-groupe $T$ qui est cyclique d'ordre $N$, et tel que pour tout $t \in T,\; t  \neq 1$, l'\'el\'ement $1-t$ soit inversible dans $A$.}
\medskip

 Un corps qui contient une racine primitive $N$-i\`eme de l'unit\'e poss\`ede cette propri\'et\'e. Si un anneau la poss\`ede, tout sur-anneau la poss\`ede aussi.
 
 Tirons les quelques cons\'equences de $({\sf P}_{N})$  qui seront utilis\'ees plus bas.
 \medskip
 
 \noindent {\bf 11.2.1.}\; {\it Dans $A[X]$, on a}
 $$
 X^N - 1 = \prod_{t\in T} (X-t).
 $$
 \medskip
 
  Cela se voit par r\'ecurrence \`a partir de la remarque suivante : soit $T' \subset T$ une partie telle que l'on ait une d\'ecomposition dans $A[X]$ en 
  $$
  X^N - 1 = P(X).\prod_{t'\in T'} (X-t') . 
  $$
  Alors, pour $t \in T,\; t \notin T'$, on a $P(t) = 0$, puisque les facteurs $(t-t')$ sont inversibles.
 \medskip
  
 \noindent {\bf 11.2.2.}\; {\it L'entier $N$ est inversible dans $A$.}
  \medskip
  
  D\'eriver l'\'egalit\'e de 11.2.1, et faire  $X = 1$.
  \medskip
  
 \noindent {\bf 11.2.3.}\; {\it Soit $t \in T$ un \'el\'ement $\neq 1$, et d'ordre divisant $d$. Alors $\sum_{j = 0}^{d-1} t^j = 0$.}
  \medskip
  
  Car cette somme est annul\'ee par l'\'el\'ement inversible $1-t$.
  \bigskip

\noindent {\bf 11.3.}\; Soit $A$ un anneau poss\'edant la propri\'et\'e  $({\sf P}_{N})$.
 Pour un groupe ab\'elien fini $G$, on pose
$$
G' \, = \, {\rm Hom}(G, T).
$$
Si $G$ est d'ordre $d$ divisant $N$,  $G'$  est un groupe ab\'elien de m\^eme ordre $d$ ({\sc N. Bourbaki} [1], A V.93)
\bigskip

Soit $u : A \rightarrow B$ un morphisme galoisien de rang $d$ divisant $N$, dont le groupe $G$ est \emph{ab\'elien}. 
\medskip
 
 Chaque homomorphisme $\theta \in G'$ est un cocycle, et d\'efinit donc le sous-$A$-module de $B$
 $$
L_{\theta}  = \{b\in B, \forall g\in G \quad   \theta(g).g(b) = b\} .
$$
Il est inversible d'apr\`es (10.5).

On v\'erifie imm\'ediatement les relations
 $$
 L_{1} = A, \qquad {\rm et}\qquad  L_{\theta}.L_{\theta'} \subset L_{\theta\theta'}.
 $$
 Ces relations permettent de munir la somme directe $\bigoplus_{\theta \in G' } L_{\theta}$ d'une structure de $A$-alg\`ebre pour laquelle l'application 
 $$
 f : \bigoplus_{\theta \in G' } L_{\theta} \; \longrightarrow \; B
 $$
 est un morphisme de $A$-alg\`ebres.
 \bigskip
 
{\bf Th\'eor\`eme 11.4.}\, {\it Soit $A$ un anneau poss\'edant la propri\'et\'e  $({\sf P}_{N})$, et soit $u : A \rightarrow B$ un morphisme galoisien de rang $d$ divisant $N$, dont le groupe $G$ est \emph{ab\'elien}. Alors, le morphisme  
$$f : \bigoplus_{\theta \in G' } L_{\theta} \; \longrightarrow \; B$$  est un isomorphisme.}
\medskip

Comme $B$ est fid\`element plat sur $A$, il suffit, d'apr\`es (9.1), de montrer que le morphisme 
$$1\otimes f : B\otimes_{A}(\bigoplus_{\theta \in G' } L_{\theta}) \; \longrightarrow \; B\otimes_{A}B$$
  est un isomorphisme. 
  
 Notons $\iota_{\theta} : L_{\theta} \longrightarrow B$ l'injection canonique, de sorte que l'on a $f(\sum_{\theta} x_{\theta}) = \sum_{\theta} \iota_{\theta}(x_{\theta})$, ce qu'on peut condenser en $f = \sum \iota_{\theta}$. Utilisons les carr\'es commutatifs introduits en 10.5.1:

$$
\begin{CD}
B\otimes_{A} L_{\theta} @>{1\otimes \iota_{\theta}}>> B\otimes_{A}B \\
 @V\omega_{\theta} VV @VV\rho V &\\
B @>>> \prod_{G}B \\
\end{CD}
$$
Passant \`a la somme, et en posant $\omega = \bigoplus \omega_{\theta}$, on obtient le carr\'e commutatif suivant o\`u $\omega$  et $\rho$  sont des isomorphismes
$$
\begin{CD}
B\otimes_{A}(\bigoplus L_{\theta}) @>{1\otimes f}>> B\otimes_{A}B \\
 @V\omega VV @VV\rho V &\\
\bigoplus_{\theta}B @>>> \prod_{G}B \\
\end{CD}
$$
Il faut voir le $B$-module $\bigoplus_{\theta}B$ comme la $B$-alg\`ebre $B[G']$ du groupe $G'$ \`a coefficients dans $B$ ; autrement dit, la composante d'indice $\theta$ du produit de $\sum x_{\theta}$ par  $\sum y_{\theta}$ est $\sum_{\theta'\theta'' = \theta} x_{\theta'}y_{\theta''}$. Il est alors clair que $\omega$ est un isomorphisme de $B$-\emph{alg\`ebres}. L'application  horizontale du bas envoie  $ \sum x_{\theta} \in \bigoplus_{\theta}B$ sur l\'el\'ement $(g\mapsto \sum_{\theta} \theta(g)^{-1}x_{\theta})$ ; en \'ecrivant la chose, on constate que c'est aussi un morphisme de $B$-alg\`ebres ; il faut montrer que c'est un isomorphisme ; le groupe $G'$ \'etant commutatif, l'application $\theta \mapsto \theta^{-1}$ est un automorphisme. Bref, on est ramen\'e \`a d\'emontrer le
\bigskip

{\bf Corollaire 11.5.}\; {\it Soit $C$ un anneau poss\`edant la propri\'et\'e  $({\sf P}_{N})$. Soit $G$ un groupe ab\'elien d'ordre $d$ divisant $N$ ; posons  $G' = {\rm Hom}(G, T)$. Alors le morphisme}
$$
f : C[G'] \; \longrightarrow \; \prod_{G} C, \qquad \theta \longmapsto (g \mapsto \theta(g))
$$
{\it est un isomorphisme.}
\medskip

Traitons d'abord le cas o\`u $G$ est cyclique (d'ordre $d$).  Le choix d'un g\'en\'erateur $g$ de $G$ et le choix d'un \'el\'ement $\zeta \in T$, d'ordre $d$, d\'eterminent un g\'en\'erateur $\theta$ de $G'$, celui d\'efini par $\theta(g) = \zeta$ ; on a donc $\theta^{i}(g^j) = \zeta^{ij}$. Les \'el\'ements $1 = \theta^0, \theta, \ldots , \theta^{d-1}$ forment une base de $C[G']$ comme $C$-module, et les applications $\delta_{j} : G \rightarrow C$, d\'efinies par $\delta_{j}(g^k) = \delta_{jk}$ forment une base de  $ \prod_{G} C$. La matrice de $f$ relativement \`a ces bases est la matrice de Van der Monde $(\zeta^{ij})$, dont le d\'eterminant $V(1, \zeta, \ldots , \zeta^{d-1}) = \prod_{i<j}(\zeta^{j}-\zeta^i) $ est inversible d'apr\`es l'hypoth\`ese $({\sf P}_{N})$.
\medskip

La d\'emonstration du cas g\'en\'eral utilise une d\'ecomposition de $G$ en produit  de groupes cycliques pour se ramener au cas pr\'ec\'edent. Si $G = G_{0}\times G_{1}$, alors $G' = {\rm Hom}(G_{0}\times G_{1}, T)$ est somme directe de ses sous-groupes $G_{0}'$ et  $G_{1}'$, et l'application \'evidente
$$
C[G_{0}']\otimes_{C}C[G_{1}'] \; \longrightarrow \; C[G']
$$
est un isomorphisme ; par suite, le morphisme $f$ de l'\'enonc\'e se factorise en
$$
C[G'] \; \simeq \; C[G_{0}']\otimes_{C}C[G_{1}']\; \stackrel{f_{0}\otimes 1}{\longrightarrow}\; \prod_{G_{0}}C\,\otimes_{C}C[G_{1}'] \;\simeq \prod_{G_{0}}C[G_{1}']  \stackrel{1\otimes f_{1}}{\longrightarrow}\; \prod_{G_{0}}(\prod_{G_{1}}C) \; \simeq \; \prod_{G}C .
$$
Cela ach\`eve la d\'emonstration du corollaire, et, par suite, celle du th\'eor\`eme.
\vspace{1cm}


\begin{center}
{\bf R\' ef\' erences}
\end{center}
\vspace{0,5cm}

\noindent [1] N. BOURBAKI, {\it Alg\`ebre}, Paris, Masson.

\noindent [2] N. BOURBAKI, {\it Alg\`ebre commutative, ch. I et II}, Paris, Masson,1961.

\noindent [3] C. GODBILLON. {\it\'El\' ements de topologie alg\' ebrique}, Paris, Hermann, 1971.

\noindent [4] E. HECKE, {\it Lectures on the Theory of Algebraic Numbers}, New-York, Springer-Verlag, 1981(Traduction de {\it Vorlesung \"uber die Theorie der algebraischen Zahlen}, 1923)

\noindent [5] M.-A. KNUS et M. OJANGUREN, {\it Th\'eorie de la descente et alg\`ebres d'Azumaya}, Springer LNM 389, 1974

\noindent [6] E. LANDAU, {\it Elementary Number Theory}, New-York, Chelsea Pub. Comp., 1958 (Traduction de {\it Elementare Zahlentheorie}, 1927)

\noindent [7] J. STILLWELL. {\it Geometry of Surfaces}, New York,
Springer-Verlag, 1992.

\noindent [8] R. SWAN. {\it Vector Bundles and Projective Modules}, Trans.
Amer. Math. Soc., Vol. 105, No 2 (Nov. 1962), 264-277.

\noindent [9] W. C.WATERHOUSE, {\it Introduction to Affine Group Schemes}, New-York, Springer-Verlag, 1979.
\vspace{2cm}

\noindent {\sc Daniel Ferrand}
\medskip

\noindent IRMAR

\noindent Universit\'e de Rennes 1,

\noindent Campus de Beaulieu

\noindent 35042 RENNES Cedex

\noindent France

\noindent e-mail : {\tt daniel.ferrand[at]univ-rennes1.fr}

\end{document}